%

%
%
\input amstex
\input amssym
%
%

\def\MediumSpaceBase{16.8truept}
\def \SingleSpaceBase{14truept}
%
%
\def \Odef{:=} 
\def \bigP{{$\Bbb P$}}
\def \Sminus{\setminus}
\def \Simage{``} 
\def \sbar#1{\bar s(#1)}
\def \sbarp#1{\bar {s'}(#1)}

\def \Tor{Todor\v{c}evi\'c}

\def \ladder#1#2#3{\langle#1_{#2} : #2 \in #3\rangle}
\def \omladder{\ladder{\varphi}{\lambda}{\Lambda}}

\def \cntord#1{#1 \in \omega_1}
\def \twocontinuum{2^{\aleph_0}}
\def \fsubl{\varphi_\lambda}
\def \symdiff{\triangle}  
\def \Inv#1{#1^{-1}}

\def \ninom{n \in \omega}
\def \forces{\Vdash}
\def \op#1#2{\langle#1, #2\rangle}  
\def \ot#1#2#3{\langle#1,\, #2,\, #3\rangle} 
\def \otp#1#2#3{\langle#1',\, #2',\, #3'\rangle} 
\def \otss#1#2#3#4{\langle#2_{#1},\, #3_{#1},\, #4_{#1}\rangle}
\def \pr#1#2{(#1,\,#2)}
\def \prxa#1{#1\pr\xi\alpha}

\def \prhb#1{#1\pr\eta\beta}

\def \prxb#1{#1\pr\xi\beta}
\def \twopre#1{{}^{#1}2}
\def \Sxh{S(\xi,\, \eta)}
\def \Sxhz{S(\xi_0,\, \eta_0)}
\def \rst{\restriction}
\def \ex{\exists}
\def \alss{\subseteq^*}  
\def \notaeq{\neq^*} 
\def \eset{\emptyset} %
%
\newcount\ThisProp\ThisProp=0
\newcount\ThisSect\ThisSect=0
\newcount\ThisDef\ThisDef=0
\newcount\ThisEqn\ThisEqn=0                              
\newcount\SectCnt\SectCnt=0
\newcount\LocalLemma
\newcount\LocalLLemma  
\newcount\LocalLLLemma
\newcount\LocalEqn
\newcount\LocalEEqn
\newcount\LocalEEEqn
\newcount\LocalEEEEqn
%
%
\def\EndPage{\vfil\eject}
\def\Section#1{\par\dimen0=\pagegoal \advance\dimen0 by -124truept
           \ifdim\pagetotal>\dimen0
                 \ifdim\pagetotal<\pagegoal  \EndPage
                 \else \vskip 8truept plus 6truept minus 4 truept \fi
           \else \vskip 8truept plus 6truept minus 4 truept \fi
           \noindent \advance\ThisSect by1
             {\bf Section \the\ThisSect  \quad} #1\par    }
\def\SubSectSkip#1{\par\vskip.1truein\noindent{\it #1}\par}
\def \Eqn {\advance\ThisEqn by1}
%
%
\def \Thm{\penalty-150\par\noindent \advance\ThisProp by1
           {\bf Theorem \the\ThisProp\quad} }
\def \endThmStmt{\par}     
\def \showThm{\penalty-150\par\noindent{\bf Proof of Theorem.\par}}
\def \endThm{\hfill $\blacksquare$  \par} 
\def \Cor{\penalty-150\par\noindent \advance\ThisProp by1
    {\bf Corollary \the\ThisProp} \quad}   
\def \showCor{\penalty-150\par\noindent{\it Proof of Corollary:\quad}}
\def \endCor{ { \hfill $\blacksquare$ } \par} 
\def \Lemma{\penalty-150\par\noindent \advance\ThisProp by1
   {\bf Lemma \the\ThisProp}\quad}
\def \showLemma{\penalty-150\par\noindent {\it Proof of Lemma:\quad}}
\def \endLemma{\hfill $\square$\par}
\def \Claim{\penalty-150\par\noindent {\it Claim:\quad}}
\def \showClaim{\penalty-150\par\noindent {\it Proof of Claim:\quad}}
\def \endClaim{ \hfill $\dashv$ \par}
\def \Prop{\penalty-150\par\noindent \advance\ThisProp by1
    {\bf Proposition \the\ThisProp}\quad}
\def \showProp{\penalty-150\par\noindent{\it Proof of Proposition:\quad}}
\def \endProp{\hfill $\square$\par} 
%
\def \DDef#1{\penalty-150\par\noindent \advance\ThisProp by1
     {\bf Definition \the\ThisProp\quad} }
\def \endDDef{\par} 
\hsize=6.5truein
\vsize=9truein
\let\dismacdefs=1  
\nopagenumbers
\baselineskip=\SingleSpaceBase
%
\font\BBX=cmbx10 scaled\magstep2
{\noindent\BBX Applications of Cohomology to Set Theory I:}

{\noindent\BBX Hausdorff Gaps}
\vskip .3 truein
{\noindent\bf Daniel E. Talayco$^*$}
\plainfootnote{}{$^*$Results in this paper consititute a substantial 
portion of my dissertation completed at The University of Michigan in 
Ann Arbor under the supervision of Andreas Blass.}

\vskip .1truein
{\noindent\it Department of Mathematics, Boston College, 
Chestnut Hill, MA  02167}

\vskip .1truein
\noindent {\bf Preprint:}  \quad Submitted to Annals of Pure and 
Applied Logic, 8/93.

\vskip .2truein
\leftskip=.3truein
\rightskip=.3truein
{\it Abstract.} We explore an application of homological algebra to set theoretic objects by
developing a cohomology theory for Hausdorff gaps.  This leads to a natural equivalence notion 
for gaps about which we answer questions by constructing many simultaneous gaps.  The first 
result is proved in ZFC while new combinatorial hypotheses generalizing $\clubsuit$ are 
introduced to prove the second result.  The cohomology theory is introduced with enough 
generality to be applicable to other questions in set theory.

Additionally, the notion of an {\it incollapsible gap} is introduced and the existence of such a
gap is shown to be independent of ZFC.  

\vskip .2 truein
\leftskip=0truein
\rightskip=0truein
\Section{Introduction}

Hausdorff gaps appear in a wide scope of applications in the literature on set theory and 
a fitting and voluminous tribute to the importance of these objects has recently been 
published by M. Scheepers, [Sc].  Hence, I shall not attempt to provide further incentive 
for their study at this time, but anticipate an interest in understanding their structure.  

In this article we explore a new approach to viewing gaps, adapting tools from cohomology 
to describe them.  In fact, although new information about gaps is gained in this way, 
this framework for investigating set theoretic structures is of as much interest to the 
author as the information about gaps.  The apparatus of cohomology will be shown to be an 
appropriate tool for understanding and directing questions about objects in set theory by 
exploring this particular example.  We will see how cohomology provides insight into the 
pertinent issues underlying gaps and from the direction so obtained will formulate and 
answer questions about gaps.  

In a future paper, we will examine the case of $\omega_1$-trees and develop a cohomology 
theory for a class of Aronszajn trees.  From this cohomological formulation of trees, we 
will be able to describe a connection between gaps and trees.  This will add evidence to 
the circumstantial similarities between gaps and trees that presently exist because of their 
constructions and similar behavior under the alchemy of forcing.  

This cohomological framework is thus initially unifying.  More enticing is the vast amount of 
technical apparatus available in the various guises of the categorical theory of homological 
algebra and derived functors.  

The modus operandi for this investigation is as follows:  1.~Describe the choices and 
interpretations of chain groups and connecting maps for the selected context; 2.~Show that 
these interpretations capture the given objects; 3.~Manipulate the cohomological apparatus 
to extract information about the associated groups; 4.~Interpret this information and use 
it to direct further investigation of the given objects.  One benefit of this approach is 
the appearance of other venues of generalization not visible from a traditionally set 
theoretic viewpoint.  Further, topological cohomology offers some useful mnemonic and 
visual associations.  

This method of exploration will be both economical and flexible.  The set theorist with 
little familiarity with cohomology should find the objects familiar enough to make the 
terminology understandable.  For the reader familiar with cohomology but unfamiliar with 
set theory, the approach should permit her or him to be introduced to the pertinent 
objects within a familiar context.  

\SubSectSkip{Notation and Conventions.}
Notation used is essentially standard for set theory.  We use [Ku] as a reference.  
Additionally, we use $\symdiff$ to denote the symmetric difference of sets.  
$P(\omega)/fin$ is considered as boolean algebra, though it should be noted that it 
maintains a group structure (actually a $\Bbb Z/2$-module structure) under the operation 
induced by $\symdiff$.  The adverb ``almost'' will mean ``modulo a finite quantity''.

The functions $\max$ and $\min$ take as arguments a list of numbers, while $\sup$, $\inf$ 
and $\sup^{+}$ (the proper supremum of a set of ordinals) take a set as an argument.  

\headline{\ifodd \pageno  \hfil {\it Gaps and Cohomology}\hfil\the\pageno
                       \else \the\pageno \hfil {\it Daniel Talayco}\hfil\fi}

$\alpha$, $\beta$, $\xi$, $\eta$, etc.,~will always denote ordinals, usually countable.  
Whenever possible, we will follow the convention that $\beta\leq\alpha$ and $\eta\leq\xi$.  
$\lambda$ will always represent a countable limit ordinal and $\Lambda$ will be the 
collection of all countable limit ordinals.  $ \kappa$ will usually denote an infinite 
cardinal.  $n$, $m$, and $i$ will usually denote elements in $\omega$.  If 
the context is clear, then quantification of these variables will implicitly be over 
objects of the proper type.  Thus in a proof where $\beta$ is used as an arbitrary 
countable ordinal, $(\forall\beta )$ will mean $ (\forall\beta\in\omega_1)$.  

\SubSectSkip{Introduction to Hausdorff Gaps.}
The usual formulation of a Hausdorff gap is as a pair of sequences of length $\omega_1$ of 
subsets of $\omega$, 
one almost increasing (increasing modulo finite) and the other almost decreasing, with the 
decreasing sequence ``above'' the increasing sequence (i.e., every element of the 
decreasing sequence almost includes every element of the increasing sequence as a subset).  
The defining property of a gap is that these sequences have no space in the middle, i.e., 
there is no set that fits between, continuing both sequences.  Hausdorff gaps, then, are a 
measure of incompleteness (in the sense of Dedekind) of $\Cal P(\omega )$.  Note that the 
order ``inclusion modulo finite sets'' on $\Cal P(\omega )$ is substantially different 
from the usual order on the reals.  Most notably here, there are uncountable well-ordered 
chains in the former.  

If each element of the top (decreasing) sequence of sets in a gap is complemented relative 
to $\omega$, the result is a pair of increasing sequences such that any pair of elements, 
one from each sequence, is almost disjoint.  The defining property then becomes:  There is 
no subset of $\omega$ continuing one sequence and remaining almost disjoint from each 
element of the other.  This is the formulation of gap used most frequently in this paper.  
Increasing sequences will be called towers.  

\DDef{$\subseteq^{*}$, $=^{*}$, etc.} Let $\subseteq^{*}$ be the pre-order on $\Cal 
P(\omega )$ given by $A\subseteq^{*}B$ if and only if $A\Sminus B$ is finite.  If this 
relationship holds, we say $ A$ is almost included in $B$ and $B$ almost includes $A$.  
If $A\subseteq^{ *}B$ and additionally $B\Sminus A$ is infinite, we say $B$ is strictly 
above $ A$ and write $A\subsetneq^{*}B$.  Denote by $A=^{*}B$ the statement that 
$A\symdiff B$ is finite and say $A$ and $B$ are almost equal.  We use this notation for 
sets in general.  

\DDef{$\noexpand \Cal P(\omega )/fin$} Let $\Cal P(\omega )/fin$ denote the Boolean 
algebra obtained by taking the subsets of the natural numbers modulo the ideal of finite 
sets.  The partial order in this algebra is induced by $\subseteq$, or equivalently, by 
$\subseteq^{*}$.  Hence we use $\subseteq^{*}$ as the order symbol on the Boolean algebra 
and may conflate elements of the Boolean algebra with representatives when convenient.  

\DDef{Almost Disjoint} Two elements of $\Cal P(\omega )$ are said to be {\it almost disjoint\/}
if their intersection is finite.  This can also be denoted $A\cap B=^{*}\emptyset$.  Unless
explicitly noted otherwise, we use the notion upwards hereditarily:  That is, two collections of
subsets of $\omega$ are said to be almost disjoint if every pair of elements, one from either
collection, is almost disjoint.

\DDef{Tower} A {\it tower}, $T=\langle T_\alpha :  \alpha\in \kappa\rangle$, in $\Cal 
P(\omega )$ is an indexed sequence of elements of $\Cal P(\omega )$ increasingly linearly 
ordered by $\subseteq^{ *}$:  $\beta <\alpha\implies T_{\beta}\subseteq^{*}T_{\alpha}$ . 
The ordinal $\kappa$ is said to be the height of the tower.  $T_{\alpha}$ will be called 
the $\alpha^{th}$ level of $T$.  A {\it subtower}, $S$, of $T$ is a tower satisfying for 
each $\alpha$, $S_{\alpha}\subseteq T_{\alpha}$ (levelwise inclusion) and additionally 
for each $\beta <\alpha$, $S_{\alpha}\cap (T_{\beta}\Sminus S_{\beta})$ is finite.  
(Restriction to lower levels of the subtower is ``faithful''.)  

The notion of a tower is more natural in $\Cal P(\omega )/fin$.  We introduce towers as 
objects in $\Cal P(\omega )$ because most of the technical work is done in this setting.  
However, it will benefit the reader to think about all statements referring to towers in 
the context of $\Cal P(\omega )/fin$.  

\DDef{Pregap} A {\it pregap\/} (or an $\omega_1$-pregap) is a pair of towers $\langle 
A,\,B\rangle$, each of height $\omega_1$ such that for each $\alpha <\omega_1$, 
$A_{\alpha}\cap B_{\alpha} =^{*}\emptyset$.  Note that both $A$ and $B$ are subtowers of 
the tower $A\cup B$, the level-wise union of $A$ and $B$.  

\DDef{Hausdorff gap} A {\it Hausdorff gap\/} is a pregap, $\langle A,\,B\rangle ,$ 
satisfying:  There is no $X$ in $\Cal P(\omega )$ such that for all $\alpha$, 
$X\supseteq^{ *}A_{\alpha}$ and $X$ is almost disjoint from $B_{\alpha}$.  We say such an 
$X$ {\it separates\/} or {\it fills\/} a pregap.  The definition is actually symmetric 
between $A$ and $B$ as can be seen by taking the complement of $X$ relative to $\omega$.  

We will usually be interested in gaps that exist inside a given tower.  Hence we give the
following definition.

\DDef{$A$ is a gap in $T$} Given a tower, $T$, in $\Cal P(\omega )$ of height $ \omega_1$, 
and a subtower of $T$, $A$, let $B$ be the levelwise complement of $A$ in $T$:  
$B_{\alpha}\Odef T_{\alpha}\Sminus A_{\alpha}$.  We say $A$ {\it is a gap in} $ T$ if and 
only if the pair $\langle A,\,B\rangle$ is a Hausdorff gap.  Assuming only that $A$ is a 
subtower of $T$, it is easy to show the tower $B$ is a subtower of $T$ and that the pair 
$\langle A,\,B\rangle$ is a pregap.  

Notice that a subtower, $A$, is a gap in a tower $T$ just in case there is no 
$X\subseteq\omega$ with $X\cap T_{\alpha}=^{*}A_{\alpha}$ for all $ \alpha$.  As above, we 
say such an $X$ fills or separates $A$.  Restricting attention to gaps inside of a tower 
does not reduce the generality of the considerations since for any Hausdorff gap $\langle 
A,\,B\rangle$, $ A$ is a gap in the tower $A\cup B$.  

The existence of a gap, proven from the axioms of ZFC, was demonstrated by Hausdorff 
in the first decade of this century.  It is this theorem which we shall generalize in 
later sections.  We produce many ``different'' gaps simultaneously, building them 
side by side.  The number of simultaneous gaps produced is eponymous, hence the $\aleph_0$ 
and $\aleph_1$ Gap Theorems.  In fact, we do substantially better in each case, at least 
for our purposes, as indicated by the corollaries following each theorem.  The motivation 
for these theorems and the notion of ``different'' are developed in the next section.  We 
first state the basic theorem and sketch the methods of generalization.  

\Thm\ (The Basic Gap Theorem) [Ha].\quad\ There is a Hausdorff gap.  \endThm

For a proof of this theorem, see [Fr, p.~36].  It is this proof from which we depart in 
the generalizations which follow.  We take a moment to review the important elements of 
the proof and indicate the directions in which changes will be made.  

Two towers are simultaneously constructed by recursion on the level.  There are two 
conflicting requirements to fulfill during the construction.  The first is a ``minimize 
intersection'' requirement to ensure that sets in the different towers are almost 
disjoint.  The second is a ``maximize intersection'' requirement to ensure the result is 
in fact a gap.  If the sets being built are $A_\alpha$ and $B_\alpha$ for 
$\alpha\in\omega_1$ then one possible maximization requirement is phrased:  $$(\forall
\alpha)(\forall r\in\omega)\,\{\beta<\alpha\,:\,A_\alpha\cap B_\beta\subseteq r\}\text{ is 
finite.}$$ This (with the minimization requirement) is sufficient to ensure the pair 
$\langle A, B\rangle$ is a gap and is used in the $\aleph_0$ Gap Theorem.  

There is an alternative to these requirements (attributed to 
\Tor)  which can be stated as follows:  $$(\forall\alpha)(\forall\beta)\,\,\,\,\,A_\alpha\cap 
B_\beta = \eset \iff \alpha=\beta.$$ This clearly implies the previous minimize 
requirement, but the proof, under this hypothesis, that the pair $\langle A, B\rangle$ is 
a gap is slightly different.  It is a variant of this second condition which is used in 
the $\aleph_1$ Gap Theorem.  

In the Basic Gap Theorem, most of the work occurs at limit stages.  At those stages, a 
first approximation to one side of the gap is defined satisfying the minimization 
requirement.  This set is carefully enlarged to satisfy the maximization requirement 
without ruining the previous work.  This is the process which is recursively applied in 
the $\aleph_0$ Gap Theorem at limit stages.  However, in the $\aleph_1$ Gap Theorem, we 
depart from this procedure adapting (the term is used loosely) the techniques of forcing 
to obtain the desired result.  

\Section{The Gap Cohomology Group}

We turn now to the motivation for the theorems which are to follow.  In fact, though 
substantial and technical, those theorems only begin to explore the potential 
generalizations and applications of homological algebra to problems in set theory and
other areas in logic.  The original observation that gaps are cohomological in nature is
due to Blass.  

The following is a list of homological ideas used but not defined here:  short and long 
exact sequences; the properties of boundary maps and how they produce cocycles and 
coboundaries; the definition of a cohomology group and how it is a measure of failed 
exactness.  The unfamiliar reader can find a reasonable introduction to this material in 
[Vi] or [H-Y].  More advanced uses of homological algebra are alluded to, but the results 
do not depend upon them.  

Herein, fix a tower, $T$, of height $\omega_1$ which may be referred to as the ambient 
tower.  Let $\Cal D\subseteq \Cal P(\omega )$ be the family of subsets of $ \omega$ 
generated by the closure under ``finite upward modification'' of sets in $T$; i.e.  
$$D\supseteq T_\alpha \land D\Sminus T_{\alpha} \text{ finite}\implies D\in \Cal D.$$ Note 
that if $\alpha ,\beta\in\omega_1$ then $T_{\alpha}\cup T_{\beta}\in \Cal D$.  More 
generally, $(\Cal D, \subseteq)$ as a partial order is directed upwards.  To each
element $D\in \Cal D$, associate the following coefficient groups:  
\item{$\bullet$}$G_D\Odef\,{}^D(\Bbb Z/2)$, the functions from $D$ into $\Bbb Z/2$, which 
will be associated with subsets of $D$, 
\item{$\bullet$}$F_D\Odef\,\bigoplus_D(\Bbb Z/2)$, the finitely supported functions from 
$D$ into $\Bbb Z/2$, associated with finite subsets of $D$, and 
\item{$\bullet$}$(G/F)_D\Odef\,G_D/F_D$, associated with the Boolean algebra $\Cal 
P(D)/fin$.  

These function groups have their group operation induced componentwise from $\Bbb Z/2$.  
We interpret such functions as characteristic functions following the convention that $f$ 
is the characteristic function of the set $f^{-1}\text{``}\{1\}$ which will be denoted by 
$\bar {f}$.  (Also, if $s(\xi )$ is a function into $\{0,1\}$ then $ \bar {s}(\xi )$ is 
the pre-image of $\{1\}$.)  In this case, the sum of two functions is the characteristic 
function of the symmetric difference of the represented sets.  

If $D_1\supseteq D_2$ then there is a natural restriction map, $\rho_ 2^1$ from $G_{D_1}$ 
to $G_{D_2}$ namely $f\,\mapsto\,f\rst D_2$; similarly for $F$ and $ G/F$.  We may now 
define cochain groups for each collection of coefficients.  A 0-cochain is just a choice 
function fixing for each $D\in \Cal D$ an element of $G_D$.  An $n$-cochain associates to 
each linearly ordered $n+ 1$ element set, $D_0\supseteq D_1\supseteq\dots\supseteq D_n$, 
an element of the ``smallest'' coefficient group, $G_{D_n}$.  Hence we define $$\roman 
C^n(\Cal D,G)\Odef\prod_{{{D_0\supseteq D_1\supseteq\dots \supseteq D_n}\atop {D_i\in \Cal 
D}}}G_{D_n}.$$ The coboundary operator connects the $n$- and $n+1$-cochain groups.  The 
elements it sends to 0 can be thought of as locally patching.  We employ functional 
notation to denote the evaluation of a cochain on its arguments.  For $x\in \roman 
C^n(\Cal D,G)$ we have $\delta x\in \roman C^{ n+1}(\Cal D,G)$ given by $$\align(\delta 
x)(D_0\supseteq&\dots\supseteq D_{n+1})\Odef\\ &\sum_{i=0}^ n\left[(-1)^ix(D_0\dots\hat 
{D_i}\dots D_{n+1})\right]+ (-1)^{n+1}\rho_{n+1}^n(x(D_0\dots D_n))\endalign$$ where $\hat 
{D_i}$ means $D_i$ is removed from the sequence.  This is essentially cohomology on a 
presheaf.  See [Jn].  The corresponding definitions for $F$ and $G/F$ are similar.  

Having made these definitions, we introduce a useful cohomological tool that relates the 
three cohomologies.  From the short exact sequence of coefficient groups, 
$$0\longrightarrow F\longrightarrow G\longrightarrow G/F\longrightarrow 0$$ we induce the 
following long exact sequence of cohomology groups:  \newcount\LESeqn \Eqn 
\LESeqn=\the\ThisEqn $$\align0\longrightarrow \roman H^0(\Cal D,F)&@>{i^{*}}>>
\roman H^0(\Cal D,G)@>{\pi^{*}}>>\roman H^0(\Cal D,G/F)@>{\delta^{*}}>>\\
&\longrightarrow \roman 
H^1(\Cal D,F)\longrightarrow \roman H^1(\Cal D ,G)\longrightarrow \roman H^1(\Cal 
D,G/F)\longrightarrow\dots\tag {\the\ThisEqn} \endalign$$ 

It is time to state a few properties of these groups and understand what is represented in the
above sequence.  Again, assume for the following that $T$ is an $\omega_1$ tower in $\Cal
P(\omega )$ and $\Cal D$ is generated from $T$ by closing under finite upward modification.  We
additionally assume that $\bigcup_{\alpha\in\omega_1}T_\alpha=\omega$, though this is
inconsequential.

\newcount\BCProp \Prop\BCProp=\the\ThisProp 
\item{(a)} $\roman H^0(\Cal D,G)\cong \Cal P(\omega )$ with $\symdiff$ (symmetric difference) 
as the group operation in $\Cal P(\omega )$.  
\item{(b)} $\roman H^0(\Cal D,F)$ depends on the structure of $T$.  If $ T$ is an inextendable 
tower then $\roman H^0(\Cal D,F)\cong fin$ (the finite subsets of $\omega$) with $\triangle$ as 
the group operation.  
\item{(c)} $\roman H^1(\Cal D,G)=0$.  

This is where the real connection between the set theoretic property of being a tower, for
example, and cohomological properties is made.  This proposition deals with the ``concrete''
case where there is no concern about modulo finite.

\showProp (a) If $x\in \roman H^0(\Cal D,G)$ and $D_0\supseteq D_1\in \Cal D$, then the 
coboundary condition on $x$ indicates $x(D_0)\rst D_1=x(D_1)$.  Thus we can define without 
ambiguity a function $\bigcup_{D\in \Cal D}x(D)$.  This is the characteristic function of 
a subset of $\omega$.  Conversely, a subset of $\omega$, $X$, induces an element of 
$H^0(\Cal D,G)$ via $x(D)=X\cap D$ for each $D\in\Cal D$ (conflating sets and their 
characteristic functions).  Since distinct subsets of $\omega$ give rise to distinct 
cocycles and every cocycle arises in this way, we have a bijective correspondence.  The 
group structure is preserved because the operation in $H^0(\Cal D,G)$ is induced by the 
operation in $\Bbb Z/2$.  

(b) If $x\in \roman H^0(\Cal D,F)$, the above applies as well.  If $T$ is an inextendable 
tower, then for any infinite $X\subseteq\omega$, there is an $\alpha$ with $T_{\alpha}\cap 
X$ infinite.  Thus the function $\bigcup_{D\in \Cal D} x(D)$ must have finite support.  
Consequently, there is always a $T$, though not necessarily of height $\omega_1$, with 
$\roman H^0(\Cal D,F)\cong fin$.  If $T$ is extendable, let $ X\subseteq\omega$ be an 
infinite set almost disjoint from all levels of $T$.  Then the cochain given by $x(D)\Odef 
D\cap X$ is a cocycle (in $\roman H^0(\Cal D,F)$).  However, distinct $ X$'s may not give 
distinct cochains.  

For one last example here, consider an extendible $\omega_1$ tower, $T$, still satisfying the
weak condition $\bigcup_\alpha T_{\alpha}=\omega$.  If $X$ is an infinite set almost disjoint
from each $T_{\alpha}$ (as $T$ is extendible) and $X'\subsetneq X$ is also infinite, then there
is an $\alpha$ with $T_{\alpha}\cap X'\neq T_{\alpha}\cap X$ (by the assumption that
$\bigcup_\alpha T_{\alpha}=\omega$).  Thus the $F$ cohomology classes associated with $X$ and
$X'$ are not equal.  As a consequence, $\vert \roman H^0( \Cal D,F)\vert =2^{\aleph_0}$; the
previous sentence demonstrates $\geq$ while the fact that an element in $\roman H^0(\Cal D,F)$
defines a subset of $\omega$ gives the reverse inequality.  It is unknown whether under
$\neg$CH other alternatives are possible.

(c) This is an example of the facility with which homological algebra can make statements 
about these structures.  The presheaf of interest, $G$, is {\it flasque\/} (also called 
flabby) which implies the higher derived functors of $\varprojlim$ (the cohomology groups 
of interest) are trivial.  See [Jn].  However, the statement can also be proved directly 
and doing so reveals how topological visualization can guide our proofs.  We wish to show 
that given a cocycle $x\in \roman C^1(\Cal D,G)$ there is an element $ y\in \roman 
C^0(\Cal D,G)$ such that $\delta y=x$.  

First, assume that $T_0\subseteq T_\alpha$ for all $\alpha$.  (Alternatively, we may assume 
$T_0=\eset$.)  We need to define $y(D)$ for 
each $D\in \Cal D$ so that\Eqn\LocalEqn=\the\ThisEqn $$(\forall D_1\supseteq 
D_2)\,\,x(D_1,\,D_2)=y(D_2)-y(D_1)\rst D_2.\tag \the\ThisEqn$$ For each $n\in\omega$ we 
will simultaneously define $y(D)(n)$ for each $D$ containing $n$.  We then verify the 
above equation holds for each $n\in D_2$.  

For each $n\in\omega$ we define $y(D)(n)$  for those $D\in\Cal D$ with $n\in D$ by:
\Eqn\LocalEEEqn=\the\ThisEqn
$$y(D)(n)\Odef -x(D,\,T_0\cup\{n\})(n).\tag \the\ThisEqn$$

Applying this definition to the right hand side of Equation (\the\LocalEqn), we have 
$$y(D_2)(n)-y(D_1)(n)=x(D_1,\,T_0\cup\{n\})(n)- x(D_2,\,T_0\cup\{n\})(n).$$ which, by the 
cochain condition on $x$ applied to the triple $D_1\supseteq D_2\supseteq T_0\cup\{n\}$, 
is equal to $x(D_1\, ,D_2)(n)$ as desired.  

If we do not assume that $T_0\subseteq T_\alpha$ for all $\alpha$, then the proof can be 
completed by defining the operation $D^*\Odef D\cup T_0$ and defining 
\Eqn\LocalEEEqn=\the\ThisEqn $$y(D)(n)\Odef (x(D^*,\,D) - x(D^*,\,T_0\cup \{n\}))(n).\tag 
\the\ThisEqn$$ The motivation for these definitions is in fact geometrical as will be 
described below.  

Suppose $D_1$ and $D_2$ are arbitrary elements of $\Cal D$ and that $n\in D_2$.  Consider 
the following diagram 
$$\CD
D_1@<{\subseteq}<<D_1^*@.\\
@V{\subseteq}VV@VV{\subseteq}V@.\\
D_2@<<{\subseteq}<D_2^*@>>{\supseteq}>T_0\cup\{n\}\endCD
$$
 The arrows represent the restriction maps induced on the
coefficient groups.  The diagram
resembles a  simplicial complex.  We
can ``push off'' the information from $x(D_1,\,D_2)$ to the other edges of
the 2-simplex containing that edge because $x$ satisfies the cochain
condition.  This was the motivation for the definition of $y$.
With this in mind, we
present the formal argument, guided by this process of ``pushing
off'' via the cochain condition.
$$\align 
(y(D_2)-y(D_1))(n)&=(x(D_2^{\ast},D_2)-x(D_2^{\ast},T_0\cup\{n\}))(n)-\\ & 
\qquad\qquad(x(D_1^{\ast} ,D_1)-x(D_1^*,T_0\cup \{n\}))(n)\\
&=x(D_2^{\ast},D_2)(n)-x(D_1^{\ast},D_1)(n)+\\ & \qquad\qquad(x(D_1^{\ast},T_0\cup
\{n\})-x(D_2^{\ast},T_0\cup\{n\}))(n)\\
&=x(D_2^{\ast},D_2)(n)-(x(D_1^{\ast},D_2)-x(D_1,D_2))(n)+x(D_1^{\ast},D_
2^{\ast})(n)\\
&=x(D_1,D_2)(n)\endalign$$
where the first equality is by definition of $y$, the second by re-arrangement, the third 
by two applications of the cochain condition on $x$ to the triples $(D_1^*,D_1,D_2)$ and 
$(D_1^*,D_2^*,T_0\cup\{n\})$, and the last equality by rearrangement and another 
application of the cochain condition on $x$ to the triple $(D_1^*,D_2^*,D_2)$.  This gives 
the desired equality and shows $y$ is sent to $x$ by $\delta$.  \endProp 

For the following proposition, we shift to the more natural context of $\Cal P(\omega )/fin$.
Recall that the partial order in $ \Cal P(\omega )/fin$ induced by $\subseteq$ is in fact
$\subseteq^{*}$ and that towers and gaps transfer to this context.  Notice that in examining
elements of $\roman H^0(\Cal D ,G/F)$, we need only consider the original tower as all other
elements of $\Cal D$ are finite modifications of elements of the tower and such finite
modification is ``washed away'' by our definitions.  

\DDef{Representative (of a cohomology class)} If $x\in \roman H^0(\Cal D,G/F)$, then $x$ is an 
equivalence class (modulo boundaries) of equivalence classes (modulo finite).  By a {\it 
representative} of 
$x(T_{\alpha} )$ I mean a subset of $\omega$ whose characteristic function is 
a representative (under the modulo finite equivalence relation) 
of the evaluation at $T_\alpha$ of a cochain representative of $x$.  

\Prop There is a bijective correspondence between $\roman H^0(\Cal D ,G/F)$ and subtowers 
of $T$ (under the equivalence relation of levelwise almost equality, i.e.  in $\Cal 
P(\omega)/fin$).  Further, if $x$ and $A$ are so related, then $A$ is a gap in $T$ if and 
only if $x$ does not lift under $\pi^*$ in the long exact sequence, (\the\LESeqn).  

Through this proposition, we see that the cohomology structure completely captures the notion
of Hausdorff gap within a given tower.

\showProp If $x\in \roman H^0(\Cal D,G/F)$, the cocycle condition on $ x$ implies that the
choices it makes (as a cochain) for subsets of each element of the tower patch in $\Cal
P(\omega )/fin$.  That is, if $A_{\alpha}$ is a representative of $x(T_{\alpha})$ for each
$\alpha$, then $\langle A_{\alpha}\,:\,\alpha\in\omega_1\rangle$ is a subtower of $T$.
Conversely, if $A$ is a subtower of $T$, define $ x(T_{\alpha})$ to be (the class containing)
$A_{\alpha}$.  These are inverse operations in $\Cal P(\omega )/fin$.

To demonstrate the second sentence of the proposition, we show the contrapositive of each
direction.  Recall that by Proposition \the\BCProp, there is a bijection between
$\Cal P(\omega )$ and $ \roman H^0(\Cal D,G)$.  Suppose $y\in \roman H^0(\Cal D,G)$ and
$\pi^{*}y=x$ where $\pi^{ *}$ is the map $\roman H^0(\Cal D,G)\to \roman H^0(\Cal D,G/F)$.  Let
$Y\subseteq \omega$ be the set given for $y$ by Proposition \the\BCProp\ and $A$
the subtower of $T$ associated with $x$ given by the preceding paragraph.  We wish to see that
$Y$ fills the pregap induced by $A$, that is, $Y\cap T_{\alpha}=^{*}A_{\alpha}$.  But this
follows immediately from the fact that $Y\cap T_{\alpha}=y(T_{\alpha})$ and $ \pi^* y=x$.

Conversely, suppose $A$ is a subtower of $T$ and $x\in \roman H^0 (\Cal D,G/F)$ is the cocycle
associated with $A$.  Suppose that $Y\subseteq\omega$ fills $ A$.  Then the cocycle $y\in
\roman H^0(\Cal D,F)$ associated with $Y$ has $ (\pi^{*}y)(T_{\alpha})$ induced by $Y\cap
T_{\alpha}$ which is almost equal to $A_{\alpha}$.  Thus $ \pi^{*}y=x$ as desired.  \endProp

From this proposition, we also see a new equivalence relation arising for gaps within a tower.
Hence, we give the following definition.

\DDef{Cohomologous gaps} We say that two gaps $A$ and $A'$ in a tower $T$ are {\it
cohomologous\/} if the levelwise symmetric difference is not a gap in $T$; i.e., if the
subtower given by $\langle A_{\alpha}\symdiff A_{\alpha}'\,:\,\alpha\in\omega_1\rangle$ is not
a gap in $T$.

By the previous proposition, this is equivalent to saying that the difference of the $G/F$
cohomology classes associated with the gaps does not lift under $\pi^{*}$.  As one example of
the reasonableness of such a relation, we have 

\Prop If $A$ and $A'$ are two gaps in $T$ which
are cohomologous and $\Bbb P$ is a notion of forcing such that it is forced by $\Bbb P$ that
$A$ is not a gap, then it is similarly forced by $\Bbb P$ that $A'$ is not a gap.
\showProp We argue in the generic extension under the stated assumptions.  Let $X$ be a set 
that fills $A$.  Thus $X\cap T_\alpha=A_\alpha$ for each $\alpha$.  Let $Y$ be a set that fills 
$A\symdiff A'$ by the assumption that $A$ and $A'$ are cohomologous.  It is straight\-forward 
to check that $X\symdiff Y$ fills $A'$.  \endProp 

As a corollary to this proposition, we have:  If $A$ is indestructible under notions of forcing
that preserve $\omega_1$ and $A$ is cohomologous to $A'$ then $A'$ is similarly indestructible.

Up to this point, we have gaps associated to the quotient of a cohomology group rather than
simply to a cohomology group.  We remedy this with the following 

\Thm The group $\roman H^1(\Cal D,F)$ is isomorphic to the set of gaps in $T$ modulo the
equivalence relation ``cohomologous'', with group operation being levelwise symmetric 
difference.  

\showThm Examining Sequence (\the\LESeqn) and Proposition
\the\BCProp\ (a) and (c) we see $$\roman H^1(\Cal D,F)\cong \roman H^0(\Cal 
D,G/F)/\ker(\delta^{*} )\cong \roman H^0(\Cal D,G/F)/\roman i\roman m(\pi^{*})$$ where the 
first $\cong$ is true because $\delta^{*}$ is a surjection (by Proposition 
\the\BCProp\ (c)) while the second is true as Sequence (\the\LESeqn) is exact.  
But the previous proposition gives the desired correspondence between this last group and the 
gaps in $T$.  \endThm

Thus, there is a cohomology group which represents gaps.  Notice that the characteristic 
properties of gaps are captured by finite sets.  As a result of this proposition, we give 
\DDef{The Gap Cohomology Group} $\roman H^1(\Cal D,F)$ 
will be called the {\it gap cohomology group}.  

Next, we state a few more properties about this equivalence relationship on gaps.

\newcount\LProp \Prop\LProp=\the\ThisProp 
\item{(a)} If $A$ is a gap in $T$ and $A'$ is (almost) obtained from $ A$ by symmetric 
difference by a constant set (i.e., there is an $X\subseteq\omega$ such that for each $\alpha$, 
$A_{\alpha}'=^{*} A_{\alpha}\symdiff(X\cap T_{\alpha})$) then $A$ and $A'$ are cohomologous.  
\item{(b)} If $A$ and $A'$ are two gaps in $T$ such that for all $\alpha$ we have 
$A_{\alpha}=^{*}A_{\alpha}'$ then $A$ and $A'$ are cohomologous.  
\item{(c)} If $A$ and $A'$ are two gaps in $T$ such that for cofinally many $\alpha$ we have 
$A_{\alpha}=^{*}A_{\alpha}'$ then $A$ and $ A'$ are cohomologous.  

\showProp For (a), we have $A_{\alpha}'\symdiff A_{\alpha}=^{*}X\cap T_{\alpha}$ for each $ 
\alpha$.  Thus $X$ fills the pregap $A'\symdiff A$.  

Part (b) follows from (a) with $X=\emptyset$.  The supposition in (c) implies the stronger
condition used in (b):  Let $\beta$ be given and let $\alpha \geq\beta$ satisfy
$A_{\alpha}=^{*}A'_{\alpha}$.  Now for each subtower $A$ and $ A'$, we know $A_{\alpha}\cap
T_{\beta}=^{*}A_{\beta}$, etc.  Thus $A_{\beta}=^{ *}A_{\alpha}\cap
T_{\beta}=^{*}A_{\alpha}'\cap T_{\beta}=^{*}A_{\beta}$.  Since $\beta$ was arbitrary, the
condition in (b) is satisfied.

\endProp 

\Section{The $\aleph_0$ Gap Theorem}

Having shown that cohomology induces an equivalence relation on the gaps within a tower, being
cohomologous, we can ask what properties the equivalence classes have.  We have seen in
Proposition \the\LProp\ that this equivalence relation smooths out some unimportant
differences in gaps.  However, it is conceivable that every pair of gaps is cohomologous.  The
main theorems of this paper, the $\aleph_0$ and $\aleph_1$ Gap Theorems indicate that this is
not the case as is explained in the corollaries following each.  The constructions, however,
are of interest in their own right, and the additional hypotheses that have arisen in
consideration of these questions seem important as well.

We now present the $\aleph_0$ Gap Theorem, so called because it is based on the construction of 
$\aleph_0$ simultaneous subtowers in a given tower.  It implies that the size of the gap 
cohomology group is at least $2^{\aleph_0}$.  

\Thm {\bf The} $\aleph_0$ {\bf Gap Theorem.  \quad }Let $T=\langle 
T_{\alpha}\,:\,\alpha\in\omega_1\rangle$ be a tower in $\Cal P(\omega )$.  Then there is an 
$\omega$ by $\omega_1$ matrix, $\langle A(m,\,\alpha )\,:\,m\in\omega ,\alpha 
\in\omega_1\rangle$, with the following properties:  

\noindent($\aleph_0$ 1) For each $m\in\omega$, $\langle A(m,\alpha
)\,:\,\alpha\in\omega_1\rangle$ is a subtower of $T$.

\noindent($\aleph_0$ 2) $(\forall m_0)(\forall m_1>m_0)(\forall\alpha )(\forall r\in\omega
)\,\{\beta <\alpha :A(m_0,\alpha )\cap A(m_1 ,\beta )\subseteq r\}$ is finite.

\noindent($\aleph_0$ 3) $(\forall\alpha\in\omega_1)(\forall\beta <\alpha )(\exists m_0)(\forall
m>m_0)\,A\pr m\alpha\cap T_{\beta} =A(m,\beta )$.

\noindent($\aleph_0$ 4) $(\forall\cntord{\alpha})\,T_{\alpha}=\bigcup \{A\pr m\alpha
:m\in\omega \}$, a disjoint union.  \endThmStmt 

The intuition behind the construction is to see $A$(column, row) with (0, 0) at the lower left.
Then each column is a tower growing upwards and each row is a partition of the associated level
of the tower $T$.  The first and second conditions ensure that each pair of columns is a
Hausdorff gap.  See [Fr, p.36] for a proof of this.  Hence each column 
is a gap in $T$.  

As previously discussed, to be a gap requires ``interaction'' (non-empty intersection) between
(sets in different) columns, a property which is assured by condition ($\aleph_0$~2).  But to
obtain the following corollary where unions of collections of columns are gaps (in particular,
are towers) requires the intersections to be controlled as formalized by condition
($\aleph_0$~3) which can be seen as a more complicated ``minimize intersection''
requirement.  This increases the delicateness and technicality required in the proof making it
reminiscent of a priority argument.  

\Cor The cardinality of the gap cohomology group is at least $\twocontinuum$.  

\showCor The third condition ensures that $$X\subseteq\omega\implies\left\langle\bigcup \{A\pr 
m\alpha :m\in X\}:\alpha\in\omega_1\right\rangle$$ is a tower by the following argument.  Let 
$A(X,\alpha )$ be the $\alpha^{th}$ level of this sequence.  We must check that 
$$\beta\leq\alpha\implies A(X,\beta )\subseteq^{*}A(X,\alpha ).$$ 

The problem is that the finite differences, $A(m,\,\beta )\Sminus A(m,\,\alpha )$, that exist
between different levels of a given column may accumulate to an infinite quantity under the
infinite union.  Examining ($\aleph_0~3$) reveals that for a fixed $\alpha$ and $\beta$ only
finitely many columns can contribute to $A(X,\beta )\Sminus A(X,\alpha )$ which hence is
finite.

If the set $X\neq\omega ,\eset$ then it is easy to see that $A(X)$ is a gap in $T$.  Recall
that two gaps are not cohomologous iff their levelwise symmetric difference is still a gap.
But it is clear that for $X,Y\subseteq\omega$ we have $A(X,\alpha )\symdiff A(Y,\alpha
)=A(X\symdiff Y,\alpha )$ and hence if $X\neq Y$ and $X\neq\omega\Sminus Y$ then $A(X)$ is not 
cohomologous with $A(Y)$.
Thus the cardinality of the gap cohomology group is at least that of the continuum.  \endCor

\showThm We construct $A\pr m\alpha$ by induction on $\alpha$, and then, for limit stages only,
induction on $m$.  I shall refer to the restrictions of ($\aleph_0$~1--4) to an ordinal
$\gamma$ (replacing $\omega_1$ by $ \gamma$ or by $\lambda$ during a limit stage) as {\it the
induction hypotheses}.  The stage $\alpha =0$ is inconsequential so long as condition
($\aleph_0$~4) is fulfilled.


\noindent {\it Successor Stage:\/} Suppose $A(m,\beta )$ have been constructed
satisfying the induction hypotheses for $\beta\leq\alpha$ and $m\in \omega$.  As $T_{\alpha
+1}\Sminus T_{\alpha}$ is infinite, it can be partitioned into infinitely many disjoint sets,
$S(m)$, $m\in\omega$.  Define $$A(m,\,\alpha +1)\Odef(A(m,\,\alpha )\cap T_{\alpha +1})\cup S(m
).$$

We check the induction hypotheses are maintained.  ($\aleph_0$~1) and ($\aleph_0$~4) are
immediate.

For ($\aleph_0$~2), fix $m_0$ and $m_1$.  We need to show that the set $\{\beta <\alpha
+1:A(m_0,\alpha +1)\cap A(m_1,\beta )\subseteq r \}$ is finite.  This set is $$\align&\subseteq
\{\beta <\alpha :(A(m_0,\alpha )\cap T_{\alpha +1})\cap A(m_1,\beta )\subseteq r\}\cup \{\alpha
\}\\ &\subseteq \{\beta <\alpha :A(m_0,\alpha )\cap A(m_1,\beta )\subseteq\sup{}^{ +}(r\cup
A(m_0,\alpha )\Sminus T_{\alpha +1})\}\cup \{\alpha \},\endalign$$ which is finite by induction
hypothesis.

For ($\aleph_0$~3) we need to show \Eqn\LocalEqn=\the\ThisEqn $$(\forall\beta\leq\alpha
)(\exists m_0)(\forall m>m_0)\,A(m,\alpha +1)\cap T_{\beta}=A(m,\beta ).\tag \the\ThisEqn$$ For
$\beta =\alpha$, this is true by the definition of $A(m,\,\alpha +1)$ and:  
\item{(a)} $(\forall m)\,S(m)\cap T_{\alpha}=\emptyset$.  
\item{(b)} $T_{\alpha}\subseteq^{*}T_{\alpha +1}\implies (\exists m_1)(\forall 
m>m_1)\,A(m,\,\alpha )\subseteq T_{\alpha +1}$.  

Now fix $\beta <\alpha$.  Then we have:  
\item{(c)} By induction hypothesis $(\aleph_0~3)$ $$(\exists m_2)(\forall m>m_2)\,A(m,\,\alpha 
)\cap T_{\beta}=A(m,\, \beta ).$$ 
\item{(d)} As $T_{\beta}\subseteq^{*}T_{\alpha}$, $S(m)$ are disjoint, and for all $ m$, 
$S(m)\cap T_{\alpha}=\emptyset$, so $$(\exists m_3)(\forall m>m_3)\,S(m)\cap 
T_{\beta}=\emptyset .$$ 

\noindent Fix $m_1$, $m_2$ and $m_3$ as in (b), (c) and (d) respectively.  Let
$m_0=\max(m_1,m_2,m_3)$ and express $A(m,\,\alpha +1)$ as $(A(m,\, \alpha )\cap T_{\alpha
+1})\cup S(m)$ to derive (\the\LocalEqn).

\noindent {\it End of Successor Stage.  


\noindent Limit Stage:\/} Now suppose $\lambda$ is a limit ordinal and for $m\in\omega$ and $
\alpha <\lambda$, sets $A(m,\,\alpha )$ have been constructed fulfilling the induction
hypotheses.  Fix a function $f:\omega\to\lambda$ which is increasing and cofinal in $\lambda$,
with $f(0)=0$.

We now construct the sets $A(m,\,\lambda )$ by recursion on $m\in \omega$.  Assume that the sets 
$A(s,\,\lambda)$ have been constructed for $s<m$.  For notation, let
$T(m,\,\lambda )$ denote the set $T_{\lambda}\Sminus\bigcup \{A(s,\,\lambda )\,:\,s<m\}$, 
the space remaining in which to build 
$A(m,\,\lambda)$.  At stage $m$, we have as induction hypotheses the 
restrictions of $(\aleph_0$~1) and $(\aleph_ 0~2)$ and additionally:

\noindent{(IH1)} $(\forall m'\geq m)(\forall\beta <\lambda )\,T(m ,\,\lambda
)\supseteq^{*}A(m',\,\beta )$.

\noindent{(IH2)} $\langle A(s,\,\lambda )\,:\,s<m\rangle$ is a disjoint family.

Fix $m\in\omega$ and assume $A(s,\lambda )$ are defined for $s<m$, satisfying the induction
hypotheses.

For $n\in\omega$, define $K(m,n)\in\omega$ to be the minimum number satisfying the following
three conditions:

\leftskip.1truein 
\item{(K1)} $K(m,n)\supseteq A(m,\,f(n))\cap\bigcup \{A(m',f(p)):  p,m'<n\,\land\,m'\neq m\}$; 
\item{(K2)} $K(m,n)\supseteq A(m,\,f(n))\Sminus T(m,\lambda )$; 
\item{(K3)} $K(m,\,n)\supseteq\bigcup \{A(m,\,f(n))\cap T_{f(p)}\Sminus A(m,\,f(p)):p<n\}$.  

\leftskip0truein To see that $K(m,\,n)$ is finite, we examine each item individually.  For
(K1), $A(m,\,f(n))\cap A(m',\,f(p))$ is finite as $m'\neq m$ and the argument of $\bigcup$ is
finite.  For (K2), this is finite by (IH1) (with $m'=m$).  For (K3), fix $n$ and $p<n$.  Then
$A(m,\,f(n))\cap T_{f(p)}=^{*}A(m,\,f(p))$ by induction hypothesis.

$K(m,\,n)$ is the amount of $A(m,\,f(n))$ to be ``removed'' in order to satisfy $A(m,\,\lambda
)\supseteq^{*}A(m,\,f(n))$.  Further, item (K1) ensures almost disjointness between sets in
different columns is maintained.  This will ensure (IH1) is maintained.

Let $B(m)\Odef\bigcup \{A(m,\,f(n))\Sminus K(m,n):\,n\in\omega \}$.  (Though including $m$ in
this notation may seem redundant, other values of $B(m)$ will be referred to later in the
proof.)  $B(m)$ is a first approximation to $A(m,\,\lambda )$.  Notice that induction
hypotheses ($\aleph_0$ 1) and (IH1) would be satisfied if we defined $A(m,\,\lambda )$ to be
$B(m)$ since:  
\item{(1)} $B(m)\supseteq^{*}A(m,\beta )$ for all $\beta <\lambda$ because for some $n$, 
$A(m,\beta )\subseteq^{*}A(m,f(n))\subseteq^{*}B(m)$, and 
\item{(2)} For $m'\neq m$ and any $i$, we have $B(m)\cap A(m',f( i))$ finite (because once 
$n>i,m'$, then $K(m,n)\supseteq A(m,\,f(n))\cap A(m',\,f(i))$ by (K3)).  This immediately gives 
for all $\beta <\lambda$, $B(m)\cap A(m',\,\beta )$ is finite.  Then by (IH1) for $ m$, for any 
$m'\geq m+1$ and any $\beta <\lambda$, we would have $A(m',\,\beta 
)\subseteq^{*}T(m+1,\,\lambda )$ as needed for (IH1) for $ m+1$.  

For $\beta<\lambda$, let $n_{\beta}$ denote the unique value of $ n$ satisfying
$f(n)\leq\beta <f(n+1)$.  Define $$J(m)\Odef\{(m',\,\beta
)\,:\,m<m'<n_{\beta}\,\land\,A(m',\,\beta )\cap B(m)\subseteq n_{\beta}\}.$$ $J(m)$ is the set
of (indices of) sets in the matrix with which $B(m)$ does not yet have ``large'' intersection
(in the sense of ($\aleph_0$ 2)).  

\Lemma$(\forall n)\,\{\beta\,:\,n_{\beta}<n\land (\exists m')\,(m' ,\,\beta )\in J(m)\}$ is
finite.  

\showLemma Fix $n\in\omega$.  It is enough to show that for each $ m'$ with $m<m'<n$
the set $pr_{m'}J(m)\Odef\{\beta\,:\,n_{\beta}<n\land (m' ,\,\beta )\in J(m)\}$ is finite since
for each $\beta$ there are only finitely many $m'$ with $(m',\,\beta )\in J(m)$.  By the
definition of $B(m)$, $A(m,f(n+1))\Sminus K(m,n+1)\subseteq B(m)$, and so $$\align
pr_{m'}J(m)&\subseteq \{\beta <f(n+1)\,:\,A(m',\,\beta ) \cap A(m,\,f(n+1))\Sminus
K(m,\,n+1)\subseteq n\}\\ &\subseteq \{\beta <f(n+1)\,:\,A(m',\beta )\cap
A(m,\,f(n+1))\subseteq\max (n,\,K(m,\,n+1))\}\endalign$$ which is finite by induction
hypothesis.  \endLemma 

Note that for any $\beta <\lambda$, if $m'>n_{\beta}$ then $(m',\, \beta )\notin J(m)$.  This
and the previous lemma imply \newcount\LocalLLL \Lemma\LocalLLL=\the\ThisProp $(\forall\beta
<\lambda )\,J(m)\cap (\omega\times\beta )$ is finite.  \endLemma

Define the function $j^m\,:\,J(m)\to\omega$ as follows.  Suppose that $(m',\,\beta )\in J(m)$.
Then let $j^m(m',\,\beta )$ be $$\inf\left[\left(A(m',\,\beta )\cap
A(m',f(n_{\beta}))\Sminus\bigcup \{T_{f(l)}:l<n_{\beta}\}\right)\cap T(m,\lambda )\Sminus
n_{\beta}\right ]$$ Note that, as $m'>m$, $A(m',\beta )\subseteq^{*}T(m,\lambda )$ by the
induction hypotheses.  (Also note that $A(m',\beta )\cap
A(m',f(n_{\beta}))=^{*}A(m',f(n_{\beta}))$ and is thus infinite.)  It is not difficult to check
in addition that the argument of $\inf$ is an infinite set.  Finally, note that $j^m(m',\,\beta
)\geq n_{ \beta}$.

Let $A(m,\,\lambda )$$\Odef B(m)\cup ran(j^m)$.  We check that the induction hypotheses are
maintained.

$(\aleph_0~1)$.  As $A\pr m\lambda\supseteq B(m)\supseteq^{*}A(m, \beta )$ for $\beta
<\lambda$, so $(\forall\beta <\lambda )\,A(m,\,\lambda )\supseteq^{*}A(m,\,\beta )$ as required
to continue the tower.  

$(\aleph_0~2)$.  Let $m'>m$ and $r\in\omega$.  We need to show that $\{\beta <\lambda :A\pr
m\lambda\cap A(m',\beta )\subseteq r\}$ is finite.  Defining $$S_n\Odef\{\beta
:n_{\beta}=n\,\land\,A(m,\,\lambda )\cap A(m',\beta )\subseteq r\},$$ we have the set of
interest equal to $\bigcup \{S_n:n\in\omega \}$. 

\Claim For all $n$, $S_n$ is finite.  

\showClaim $$\align S_n&\subseteq \{\beta <f(n+1):A\pr m\lambda\cap A(m',\beta )\subseteq r\}\\ 
&\subseteq \{\beta <f(n+1):B(m)\cap A(m',\beta )\subseteq r\}\\ &\subseteq \{\beta 
<f(n+1):A(m,f(n+1))\Sminus K(m,n+1)\cap A(m',\beta )\subseteq r\}\\ &\subseteq \{\beta 
<f(n+1):A(m,f(n+1))\cap A(m',\beta )\subseteq\max (r,K(m,n+1))\}\\ \endalign$$ which is finite 
by induction hypothesis.  \endClaim 

\Claim If $n>\max(r,m')$ then $S_n=\emptyset$.  

\showClaim Suppose $\beta\in S_n$ where $n>\max(r,m')$.  Then $n_{\beta}=n$ and $A\pr 
m\lambda\cap A\pr{m'}\beta\subseteq r\subseteq n$.  Since $B(m)\subseteq A\pr m\lambda$ this 
gives $B(m)\cap A\pr{m'}\beta\subseteq n$.  But $m<m'<n$, and so $(m',\,\beta )\in J(m)$.  
Consequently, $j^m(m',\,\beta )\in A\pr m\lambda\cap A(m',\beta )\Sminus n$, which contradicts 
the deduction that this set is empty.  \endClaim 

Thus ($\aleph_0~2$) is verified.  

(IH1).  Fix $m'>m+1$ and $\beta <\lambda$.  We wish to show that $T(m+1,\,\lambda
)\supseteq^{*}A(m',\,\beta )$.  Since by induction hypothesis, $T(m,\,\lambda
)\supseteq^{*}A(m',\,\beta )$, it is sufficient to show that $ A(m,\,\lambda )\cap A(m',\,\beta
)$ is finite.  Further, we may assume $\beta$ is of the form $f(n)$ for some $n\in\omega$.
(See the proof of ($\aleph_ 0~3$) below for the proof of a similar
statement.)

Assume $\beta =f(n)$.  It was noted that $B(m)\cap A(m',\beta )$ is finite.  To establish
$ran(j^m)\cap A(m',\beta )$ is finite, we show $$ran(j^m)\cap A(m',\beta )\subseteq
ran(j^m\restriction[J(m)\cap (\omega\times f(n_{\beta}+1))])$$ and appeal to Lemma
\the\LocalLLL\ above.  Let $(p,\,\gamma )\in J(m)\Sminus(\omega\times
f(n_{\beta}+1))$.  Then $\gamma\geq f(n_{\beta}+1)$ and so $n_{\gamma}>n_{\beta}$.  By the
definition of $j^m$, we have $$j^m(p,\,\gamma )\notin T_{f(n_{\beta})}=T_{\beta}\supseteq A(m'
,\,\beta )$$ as desired.  This completes the proof for (IH1).

(IH2) is immediate since $A(m,\,\lambda )\subseteq T(m,\,\lambda )$ which is disjoint from
$A(s,\,\lambda )$ for $s<m$.

\noindent {\it End of Limit Construction.}

It remains to check that $(\aleph_0~3)$ and $(\aleph_0~4)$ are satisfied after the completion
of the construction of $A(m,\,\lambda )$ for $m\in\omega$.

($\aleph_0~3$).  We must show that $$(\forall\beta <\lambda )(\exists m_0)(\forall
m>m_0)\,A(m,\,\lambda )\cap T_{\beta}=A(m,\beta ).$$

\Claim It is sufficient to check this for $\beta$ of the form $f(n)$.  For suppose it holds of
such ordinals and $\beta\in\omega_1$ is arbitrary.  Fix $ n$ such that $\beta <f(n)$.  Then
there are $m_0$, $m_1$ and $m_2$ such that:  
\item{(1)} $(\forall m>m_0)\,A\pr m\lambda\cap T_{f(n)}=A(m,\,f(n))$, i.e., ($\aleph_0~3$) 
holds for $f(n)$, 
\item{(2)} $(\forall m>m_1)\,A(m,\,f(n))\cap T_{\beta}=A(m,\beta )$, which holds by induction 
hypothesis, 
\item{(3)} $(\forall m>m_2)\,A\pr m\lambda\cap (T_{\beta}\Sminus T_{f (n)})=\eset$ which is 
possible as $T_{\beta}\Sminus T_{f(n)}$ is finite and the $A\pr m\lambda$ are disjoint.  

Then for $m>\max(m_0,m_1,m_2)$, $$\align A\pr m\lambda\cap T_{\beta}&=A\pr m\lambda\cap
[(T_{\beta}\Sminus T_{f(n)})\cup (T_{\beta}\cap T_{f(n)})]\\ &=\eset\cup (A\pr m\lambda\cap
T_{f(n)}\cap T_{\beta})\text{ as } m>m_2\\ &=A(m,\,f(n))\cap T_{\beta}\text{ as }m>m_0\\
&=A(m,\beta )\text{ as }m>m_1.\endalign$$ \endClaim 

We use induction on $n$.  Fix $n\in\omega$ and assume the claim holds for all $p<n$.  We wish to 
show \Eqn $$(\ex m_0)(\forall m>m_0)\,A\pr m\lambda\cap T_{f(n)}=A(m,\,f(n)).\tag \the\ThisEqn$$ 
Since there are only finitely many $p<n$, we have by induction hypothesis that\Eqn 
\newcount\EQOne\newcount\EQTwo\EQOne=\the\ThisEqn $$(\exists m_1)(\forall p<n)(\forall 
m>m_1)\,A(m,\,f(n))\cap T_{f (p)}=A(m,\,f(p)).\tag \the\ThisEqn$$ Fixing such an $m_1$, this 
gives \Eqn\EQTwo=\the\ThisEqn $$(\forall p\leq n)(\forall m>m_1)(\forall m'\neq m)\,A(m,\,f(n)) 
\cap A(m',\,f(p))=\emptyset\tag \the\ThisEqn$$ (where we use additionally that $A(m,\,f(p))\cap 
A(m',\,f(p))=\eset$).  Next we have $$(\exists m_2)(\forall m>m_2)A(m,\,f(n))\subseteq 
T_{\lambda}.$$ Let $m_0\Odef\max(m_1,m_2,n)$.  

\Claim$(\forall m>m_0)\,A\pr m\lambda\cap T_{f(n)}=A(m,\,f(n))$.  
\showClaim This is demonstrated by establishing the following three facts.  
\item{(A)} $(\forall m>m_0)\,\,ran(j^m)\cap T_{f(n)}=\eset$.  
\item{(B)} $(\forall m>m_0)\,\,B(m)\cap T_{f(n)}\subseteq A(m,\,f(n) )$.  
\item{(C)} $(\forall m>m_0)\,\,K(m,n)=0.$ 

(C) implies $(\forall m>m_0)\,B(m)\supseteq A(m,\,f(n))$ which combined with (B) gives
$(\forall m>m_0)\,B(m)\cap T_{f(n)}=A(m,\,f(n))$.  With (A) we can deduce the claim since $A\pr
m\lambda =B(m)\cup ran(j^m)$.

\noindent {\it Proof of\/} (A).\quad Show $(\forall m>m_0)\,\,ran(j^m)\cap T_{f(n)}=\eset$.

Let $(m',\,\beta )\in J(m)$.  We may infer that $n\leq m_0<m<m'<n_{ \beta}$, where the first
inequality holds by the definition of $m_0$ and so $f(n)<f(n_{\beta})\leq\beta$.  Thus
$j^m(m',\,\beta )\notin T_{f(n)}$ by the definition of $j^m$.  This establishes (A).

\noindent {\it Proof of\/} (B).\quad Show $(\forall m>m_0)\,\,B(m)\cap T_{f(n)}\subseteq
A(m,\,f(n))$ or, from the definition of $B(m)$, show that
$$\bigcup_{p'\in\omega}\left(A(m,\,f(p'))\cap T_{f(n)}\Sminus K(m ,p')\right)\subseteq
A(m,\,f(n)).$$ We examine the two cases when $p'\leq n$ and when $p'>n$.  In the first case we
appeal to display (\the\EQOne) which implies $(\forall m>m_0) (\forall p\leq
n)\,A(m,f(p))\subseteq A(m,\,f(n))$.

When $p'>n$, $K(m,\,p')\supseteq [A(m,\,f(p'))\cap T_{f(n)}\Sminus A(m,\,f(n))]$ by the third
part of the definition of $K(m,\,p')$.  Consequently $A(m,\,f(p'))\cap T_{f(n)}\Sminus
K(m,\,p')\subseteq A(m,\,f(n))$.  Hence the displayed union over $p'$ is contained in the
desired set, $ A(m,\,f(n))$.  This establishes claim (B).

\noindent {\it Proof of\/} (C).  Show $(\forall m>m_0)\,\,K(m,n)=\eset$.

There are three parts to the definition of $K(m,n)$.  It is sufficient to show that for $m>m_0$
the right hand side of each part is empty.

For (K1) we must show $$(\forall m>m_0)\,\left(A(m,\,f(n))\cap\bigcup \{A(m',f(p)):p,m'<n\text{
and } m'\neq m\}\right)=\eset.$$ This is immediate from the definition of $m_0$ and display
(\the\EQTwo).

For (K2), show $(\forall m>m_0)\,A(m,\,f(n))\Sminus T(m,\lambda )=\eset$.  Using the definition
of $T\pr m\lambda$, and since $m>m_0$ implies that $A(m,\,f(n))\subseteq T_{\lambda}$, this
reduces to $(\forall m>m_0)(\forall s<m)\,\,A(m,\,f(n))\cap A\pr s\lambda =\eset$.  As
$A(s,\lambda )=B(s)\cup ran(j^s)$, we show 
\item{(i)} $(\forall m>m_0)(\forall s<m)\,A(m,\,f(n))\cap ran(j^s)=\eset$, and 
\item{(ii)} $(\forall m>m_0)(\forall s<m)\,A(m,\,f(n))\cap B(s)=\eset$.  

For (i), suppose for arbitrary $\beta ,\,m',$ that $(m',\,\beta )\in J(s)$.  We show that
$j^s(m',\,\beta )\notin A(m,\,f(n))$ in two cases.  First, if $n_{\beta}\leq n$ then $m\geq
m_0\geq n\geq n_{\beta}>m'$ (where the last inequality follows from the definition of $J(s)$)
gives $m\neq m'$ which by display (\the\EQTwo) gives $A(m,\,f(n))\cap
A(m',\,f(n_{\beta}))=\eset$.  But $j^s(m',\, \beta )\in A(m',\,f(n_{\beta}))$.  In the second
case, $n_{\beta}>n$.  From the definition of $j^s,\, j^s(m',\,\beta )\notin T_{f(n)}\supseteq
A(m,\,f(n))$.  This proves (i).

For (ii), note that $$A(m,\,f(n))\cap B(s)=\bigcup_{p\in\omega}\left[A(m,\,f(n))\cap
A(s,f(p))\Sminus K(s,p)\right],$$ and again we have two cases.  If $p\leq n$ then immediately
$A(m,\,f(n))\cap A(s,\,f(p))=\eset$ by display (\the\EQTwo) and the fact that $s\neq m$.

In the second case, where $p>n$, we have by the third clause in the definition of $K(s,p)$ that
$$\align K(s,p)&\supseteq A(s,f(p))\cap T_{f(n)}\Sminus A(s,f(n))\\ &\supseteq A(s,f(p))\cap
A(m,\,f(n))\Sminus A(s,f(n))\\ &\supseteq A(s,f(p))\cap A(m,\,f(n)).\endalign$$ where the last
line is true because $A(m,\,f(n))\cap A(s,\,f(n))=\eset$.  Thus $A(s,f(p))\cap
A(m,\,f(n))\Sminus K(s,p)=\eset$ as desired.  This establishes (ii) which completes the proof
for (K2).

For (K3), we show that $$\bigcup_{p<n}\left[A(m,\,f(n))\cap T_{f(p)}\Sminus A(m,f(p))\right
]=\eset.$$ \par This follows immediately from display (\the\EQOne) since
$$A(m,f(p))=A(m,\,f(n))\cap T_{f(p)}\implies A(m,f(n))\cap T_{f(p )}\Sminus A(m,f(p))=\eset.$$
The premise of this display holds when $m>m_0$.  This completes the proof for (K3) which
finishes claim (C) that $K(m,n)=\eset$ for all but finitely many $m$.  \endClaim

By the reasoning presented after the statements of (A)---(C), we conclude the claim that for
all but finitely many $m$, $A\pr m\lambda\cap T_{f(n)}=A(m,\,f(n))$.  This completes the proof
that ($\aleph_0~ 3$) holds through the limit stage.

($\aleph_0~4$).  The sets $A\pr m\lambda ,\,m\in\omega$ are disjoint by construction.  Suppose
that their union does not exhaust $T_{\lambda}$.  Notice that the quantity of $T_{\lambda}$
remaining must be almost disjoint from $T_{\beta}$ for all $\beta <\lambda$ (and hence almost
disjoint from each $A(m,\, \beta )$ for $m\in\omega$ and $\beta <\lambda$).  This follows from
($\aleph_0~ 3$).  Consequently, $A(0,\,\lambda )$ can be expanded to contain this set without
affecting the other hypotheses.  Note in particular that $\langle
A\pr0\alpha\,:\,\alpha\leq\lambda\rangle$ will be continue to be a subtower of $T$.

\noindent {\it End of Limit Stage.}

This completes the construction of the desired $\omega\times\omega_ 1$ matrix.  Since the
properties of this matrix are all stated with quantifiers over countable ordinals, the proofs
of the persistence of the induction hypotheses through the recursion establishes that the
matrix has the stated properties.  \endThm

Immediate attempts to improve this theorem were resisted by apparently combinatorial
complications.  These difficulties had the flavor of independence results and indeed the only
successful attacks on the problem have relied on combinatorial principles which are known to be
consistent with and independent of the axioms of ZFC.  At this point, it
is necessary that these principles be introduced in their proper context.

\Section{New Combinatorial Hypotheses}

We turn now to the combinatorial hypothesis that will be used to prove the $\aleph_1$ Gap
Theorem.  This and related hypotheses seem to be of interest in their own right, and so I take
the opportunity to prove some statements about their relative consistency.

An important type of object for these definitions is the following.  The reader can find further 
information on such objects in [D-S].  
\DDef{Ladder System} A {\it ladder system,} $\ladder\varphi\lambda\Gamma$, on a set of limit 
ordinals of countable cofinality, $\Gamma$, is a $\Gamma$-indexed collection of increasing 
$\omega$-sequences, $\fsubl$, each cofinal in its respective $\lambda$.  \endDDef 

Recall the definition of $\clubsuit$ from [Os]:  
\DDef{$\clubsuit$} $\clubsuit$ is the statement that there is a ladder system, $\omladder$, such 
that for every uncountable set $X\subseteq\omega_1$ there is a $\lambda\in\Lambda$ with 
$\fsubl\subseteq X$.  \endDDef 

Most of the hypotheses involved follow the basic form of $\clubsuit$.  That is, they state the
existence of a sequence of sets having some property with respect to other sets.  In general,
we will follow the convention that a sequence satisfying these properties is called a
$\clubsuit$-sequence ($\diamondsuit$-sequence, H2-sequence, etc.)

The reader is refered to [Ku, p.~80] for the statements of the hypotheses in the $\diamondsuit$
family.

Blass has pointed out that H2, defined below, can be phrased as a negative partition relation
connecting these ideas to the work of \Tor\ and others.  This seems to reflect the implicit
connection between the properties used to ensure a pair of towers is a gap---an event occurring
between sets at different levels---and partition relations on pairs of ordinals.  In addition,
it has led to weakened forms of the hypotheses, also given below, which are more easily seen
to be independent of ZFC.

\DDef{H0} H0 is the statement that there is a ladder system $\ladder\varphi\lambda\Lambda$ such
that for each stationary subset $S$ of $\omega_1$ there is a $\lambda\in S$ such that
$\fsubl\subseteq S$.

Compared to $\clubsuit$, we have strengthened the statement in requiring the ``self-reference''
of $S$, while weakening the universal quantifier to stationary sets.  In any case, this turns
out to be inconsistent with ZFC.

\Claim H0 is not consistent with ZFC.  
\showClaim Suppose $\ladder\varphi\lambda\Lambda$ were an H0-sequence.  Inductively define a 
set $S$ such that $$\lambda\in S\iff\varphi_{\lambda}\nsubseteq S$$ We show that $S$ is 
stationary, immediately contradicting H0.  Let $C$ be a club and suppose $C\cap S=\eset$.  By 
the assumption of H0, there is a $\lambda\in C$ such that $\fsubl\subseteq C$.  But then 
$\lambda\in S$, a contradiction.  Hence $S$ is stationary.  \endClaim 

Fortunately, the same fate does not befall the following weakenings of H0.

\DDef{H1} H1 is the statement that there is a ladder system $\ladder\varphi\lambda\Lambda$ such
that for each stationary subset $S$ of $\omega_1$ there is a $\lambda\in S$ such that
$\vert\fsubl\cap S\vert =\aleph_0$.  \endDDef

\DDef{H2} H2 is the statement that there is a ladder system $\ladder\varphi\lambda\Lambda$ such
that for each stationary subset $S$ of $\omega_1$ there is a $\lambda\in S$ such that
$\fsubl\cap S\neq\eset$.  \endDDef

\Prop $\diamondsuit^{*}\implies$ H1 $\implies$ H2.  In particular, H1 and H2 are consistent
with ZFC.  

\showProp\ The second implication is immediate.  Let $\langle\Cal 
D_{\alpha}:\alpha\in\Lambda\rangle$ be a $ \diamondsuit^{*}$-sequence.  For each 
$\lambda\in\Lambda$ define $\varphi_{\lambda}$ a cofinal $\omega$-sequence in $\lambda$ such 
that for each $D\in \Cal D_{\lambda}$ which is cofinal in $\lambda$, $\varphi_{\lambda}\cap D$ 
is infinite.  This is done by enumerating the $ D\in \Cal D_{\lambda}$ which are cofinal in 
$\lambda$ and recursively defining $\varphi_{\lambda}$.  If there are no cofinal elements in 
$\Cal D_{\lambda}$, then let $\varphi_{\lambda}$ be arbitrary.  

I claim that this $\varphi_{\lambda}$ sequence is an H2-sequence.  For let $S$ be a stationary 
subset of $\omega_1$ and let $C$ be a club as in the definition of the 
$\diamondsuit^{*}$-sequence, that is, where $S$ is predicted.  Let $\lambda\in C\cap S$ such 
that $ C\cap S$ is cofinal in $\lambda$.  This is possible because $C\cap S$ is stationary.  As 
$\lambda\in C$, we know $S\cap\lambda\in \Cal D_{\alpha}$ and $S\cap\lambda$ is cofinal in $ 
\lambda$.  By the definition of $\varphi_{\lambda}$, we have the desired statement that $ 
\varphi_{\lambda}\cap S$ is infinite.  This shows that H1 and H2 hold in L and are thus 
consistent with ZFC.  \endProp 

Of course, $\diamondsuit\implies\text{CH}$, and we are interested in statements about 
gaps when CH does not hold, too.  Further, the use of $\diamondsuit^{*}$ in the above proof 
seems to be more than is necessary.  It would be more satisfying to have a better understanding 
of the power of H2.  Towards this end, we will show H2 is consistent with $\neg$CH.  In fact, 
the proof below shows H2 is consistent with the continuum being anything reasonable and can be 
easily adapted to show the same for H1.  

This theorem is proved by showing that an H2-sequence is preserved under notions of forcing
that satisfy an apparent strengthening of a previously known condition.  We begin by stating
this new condition and proving lemmas that will help show familiar notions of forcing satisfy
the condition.

\DDef{Property SK} We say that a notion of forcing, $ \Bbb P$, has {\it property SK} if and
only if for any sequence of conditions $\langle p_{\alpha}\,
:\,\alpha\in S\rangle$ indexed by a stationary set $S\subseteq\omega_1$ there is a stationary
$T\subseteq S$ such that for all $\alpha$, $\beta\in T$, $p_{\alpha}$ and $p_{\beta}$ are
compatible.  \endDDef

SK can be read as strong Knaster or stationary Knaster as this is a strengthening of property
K.

\Lemma\ If $\langle A_{\alpha}\,:\,\alpha\in S\rangle$ is a collection of finite subsets of $
\omega_1$ with $S\subseteq\omega_1$ stationary, then there is a stationary $ T\subseteq S$ such
that $\langle A_{\alpha}\,:\,\alpha\in T\rangle$ is a $\Delta$-system.
\newcount\Dlemma\Dlemma=\the\ThisProp

\showLemma\ (Blass)\quad Thinning $S$, we may assume all $A_\alpha$ have the same cardinality,
$n$; for $k<n$ let $a_\alpha(k)$ be the $k^{th}$ element of $A_\alpha$.  If there is a
stationary set of $\alpha$'s for which $a_\alpha(n-1)$ is bounded, then for these $\alpha$'s there 
are only
countably many different $A_\alpha$'s and so stationarily many are the same.

Otherwise, let $k$ be the least number such that $a_\alpha(k)$ is unbounded on every stationary
set of $\alpha$'s.  Note that the same must hold for all $i$ between $k$ and $n$.  As above, we
can thin the index set to a stationary set such that $\{a_\alpha(0),\dots,a_\alpha(k-1)\}$ is
independent of $\alpha$.  As $a_\alpha(k)$ cannot be a regressive function of $\alpha$ on any
stationary set, we can thin to arrange that $\alpha\leq a_\alpha(k)$ for all $\alpha$.
Further, by intersecting the index set with a suitable club (namely, $\{\alpha\,:\,(\forall
\beta<\alpha)\, a_\beta(n-1)<\alpha\}$) we have $a_\beta(n-1)<\alpha$ whenever $\beta<\alpha$.
This collection of $A_\alpha$'s forms a $\Delta$-system with kernel
$\{a_\alpha(0),\dots,a_\alpha(k-1)\}$.  \endLemma

Now we demonstrate that two of the most familiar notions of forcing have property SK.

\Lemma\ If $\kappa\geq\omega_1$ and $\Bbb P$ is the set of finite partial functions from $
\kappa$ into 2, then $\Bbb P$ has property SK.  (I.e., Cohen forcing has property SK.)
\showLemma Let $p_{\alpha}\in \Bbb P$, $\alpha\in S$, where $S\subseteq\omega_ 1$ is stationary
be given.  Notice that the cardinality of the set of all finite partial functions on the union
of the domains of the $p_{\alpha}$ is $\aleph_1$.  As the previous lemma addresses only
properties of extensionality and cardinality, we can apply it to $\langle
dom(p_{\alpha})\,:\,\alpha\in S\rangle$ to get a $ \Delta$-system on a stationary $T'\subseteq
S$.  We may then reduce to a set $T\subseteq T'$ such that the restriction of $p_{\alpha}$ to
the root of the $\Delta$-system is independent of $\alpha$ in $T$.  Then the union of any two
conditions indexed by $T$ is a common extension of each.  \endLemma

\Lemma\ Random forcing has property SK.

The proof actually shows that any $\sigma$-linked forcing has property SK.  
\showLemma\ We use
the following fact:  If $p$ is a Borel set of positive measure, then for almost all $x\in p$,
the density of $p$ in intervals around $x$ goes to 1 as the interval width goes to 0. (This
result, known as the Lebesgue density theorem, can be established by showing sets without this
property have measure 0.)  So for $q\in \Bbb Q$, the rationals, and $n\in\omega$, consider the
set $$S^q_n\,:=\big\{p\,:\,\text{The density of }p\text{ in }\big(q -\frac 1n,q+\frac
1n\big)\text{ is greater than }\frac 12\big\}.$$ Notice that any pair of conditions in this set
are compatible (have intersection with positive measure.)  Secondly, for any $p$, there is a
$q\in \Bbb Q$ and an $n\in\omega$ with $p\in S^q_n$.  If $S$ is stationary and indexes a set of
conditions, then there must be a $q$ and an $n$ with $\{\alpha\,:\,p_{\alpha}\in S^q_n\}$
stationary.  This completes the lemma.  \endLemma

\Lemma\ Having property SK is preserved under finite support iteration.  I.e., if $\langle
Q_i,\,\pi_i\,\,:\,i\in\kappa\rangle$ is a finite support iteration and $\forces_{\Bbb
P_\gamma}\text{``}\dot Q_\gamma$ has property SK'' for all $\gamma<\kappa$, then the resulting 
partial order has
property SK.  
\showLemma Let $\langle Q_i,\,\pi_i\,\,:\,i\in\kappa\rangle$ be a finite support
iteration and $ p_{\alpha}\in \Bbb P$ for $\alpha\in S$ where $S\subseteq\omega_1$ is
stationary.  Since $spt(p_{\alpha}) \in [\kappa ]^{<\aleph_0}$, we may apply Lemma
\the\Dlemma\ to get a stationary $T_0\subseteq S$.  For each $ \gamma$ in the
root of this $\Delta$-system, successively get $T_{n+1}\subseteq T_n$ stationary such that
$$\forces_{\Bbb P_\gamma}(\forall\alpha ,\beta\in T_{n+1})\,\,\text{``}p_{\alpha} (\gamma
)\text{ is compatible with }p_{\beta}(\gamma )\text{''},$$ by applying the fact that
$Q_{\gamma},\,\,\pi_{\gamma}$ has property SK.  As the original root was a finite set, there is
a single $T_N$ demonstrating the lemma.  \endLemma

We are finally ready to state and prove the theorem.

\Thm Con(ZFC + H2 + $\neg$CH).  That is, if ZFC is consistent, then so is ZFC + H2 + $\neg$CH.  

\showThm
Start with a model of H2, for example a model of $\diamondsuit^{*}$.  Add $\aleph_ 2$ Cohen
reals.  I intend to show that the H2-sequence in the ground model continues to enjoy this
property in the extension.  

Let $\dot {S}$ be the name of a stationary subset of $\omega_1$ and fix $ p\in \Bbb P$.  Define 
$p_{\alpha}\,:=\Vert\alpha\in\dot {S}\Vert\land p$.  As $\forces\text{``}\dot { S}\text{ is 
stationary''}$, we have $\{\alpha\,:\,p_{\alpha}\text{ is not 0}\}$ is stationary in the generic 
extension (associated with any generic set containing $p$) as it contains $\dot {S}$.  Hence it 
is stationary in the ground model where it is definable.  Since $\Bbb P$ has property SK, there 
is a stationary $T\subseteq\omega_1$ such that $\{p_{\alpha}\,:\,\alpha\in T\}$ is a pairwise 
compatible set.  By H2 in the ground model, there is a $\lambda\in T$ such that 
$\varphi_{\lambda} \cap T\neq\emptyset$.  So for some $\alpha\in\varphi_{\lambda}$, we have 
$p_{\alpha}$ and $p_{\lambda}$ are compatible, $ p_{\alpha}\land p_{\lambda}\leq p$ and 
$$p_{\alpha}\land p_{\lambda}\forces\lambda\in\dot {S}\,\land\,\dot { 
S}\cap\varphi_{\lambda}\neq\emptyset .$$ Since we started with an arbitrary condition and an 
arbitrary name for a stationary set and we found an extension of the condition which forces that 
the set has non-empty intersection with some $\varphi_{ \lambda}$ and forces that $\lambda$ is 
in the stationary set, we see that the H2 sequence in the ground model remains an H2 sequence in 
the extension.  We have thus completed the proof of the theorem.  \endThm 

Notice $\aleph_2$ Cohen reals in the above proof could be replaced with any number of Cohen or
random reals (added in a finite support iteration), or any other notion of forcing known to
have property SK.

Since the development of H2 and the above discourse, Blass has noted that a weakening of H2 is
all that is really needed in the proof of the $\aleph_1$ Gap Theorem below.  Additionally, this
weakening follows from $\diamondsuit$ rather than $\diamondsuit^*$.

\DDef {weak-H1} Weak-H1 is the statement that there is a ladder system
$\ladder\varphi\lambda\Lambda$ such that for all functions $f\,:\,\omega_1\to\omega$ there is a
limit $\lambda$ such that $\{\alpha\in\varphi_\lambda\,:\,f(\lambda)=f(\alpha)\}$ is infinite.
\endDDef

\DDef {weak-H2} Weak-H2 is the statement that there is a ladder system
$\ladder\varphi\lambda\Lambda$ such that for all functions $f\,:\,\omega_1\to\omega$ there is a
limit $\lambda$ and an $\alpha\in\varphi_\lambda$ such that $f(\lambda)=f(\alpha)$.  \endDDef

Notice that H1 implies weak-H1 and H2 implies weak-H2.  These can be phrased in terms of
partitions as well.  For example, for a weak-H2 sequence, given any partition of $\omega_1$
into countably many pieces, there is necessarily one piece which contains both some $\lambda$
and an element of that $\fsubl$.

\Prop (Blass) $\diamondsuit \implies$ weak-H1 $\implies$ weak-H2.
\showProp Fix a $\diamondsuit$-sequence of functions $f_\alpha\,:\,\alpha\to\omega$.  Let 
$\fsubl$ be an $\omega$-sequence increasing to $\lambda$ additionally satisfying $$(\forall 
n\in\omega)\text{ If } f_\lambda^{-1}\{n\}\text{ is cofinal in } \lambda\text{ then it meets 
}\fsubl\text{ in an infinite set.}$$ Such a $\varphi_\lambda$ may be constructed recursively by 
infinitely often addressing each $n$ for which $f_\lambda^{-1}\{n\}$ is cofinal in $\lambda$.  
We show this sequence satisfies weak-H1.  

Let $g:\omega_1\to\omega$.  Let $A\Odef\{n\in\omega\,:\,g^{-1}\{n\}$ is uncountable$\}$.
Define $\beta\Odef\sup^+\bigcup\{g^{-1}\{n\}\,:\,n\in\omega\Sminus A\}$.  Notice that
$\beta<\omega_1$ and that if $\xi>\beta$ then $g(\xi)\in A$.

For each $n\in A$, the set of limit points of $g^{-1}\{n\}$ is a club; let $C$ be the
intersection of these countably many clubs.  Because the $f_\alpha$'s form a
$\diamondsuit$-sequence, there is a limit ordinal $\lambda\in C$ with $\lambda>\beta$ and
$g\rst\lambda=f_\lambda$.

Fix such a $\lambda$ and let $n=g(\lambda)$.  By the above remark, $n\in A$.  As $\lambda\in
C$, it follows that $\lambda$ is a limit point of $g^{-1}\{n\}\cap\lambda =
f_\lambda^{-1}\{n\}$.  So, by construction, $\fsubl$ meets $f_\lambda^{-1}\{n\}$, say at
$\gamma$.  Then $g(\gamma)=f_\lambda(\gamma)=n$ and $g(\lambda)=n$ and $\gamma\in\fsubl$ as
required.

The second implication is immediate from the definitions.  \endProp

There is another reason for the introduction of the weak forms of these hypotheses.  Given a
ladder system, it is relatively easy to describe a c.c.c.\ notion of forcing which adjoins, via
$\aleph_1$ many dense sets, a function from $\omega_1$ to $\omega$ which demonstrates that the
given ladder system is not a weak-H2 sequence.  Consequently, we have

\Prop Under MA, weak-H2 does not hold; hence neither do weak-H1, H2 nor H1 hold.  \endProp

Before leaving these hypotheses, there is one more curiousity to point out.  It is well known
that in the conclusion of $\diamondsuit$, it is equivalent to assume that the
$\diamondsuit$-sequence guesses correctly just once or to assume that the set of correct
guesses is stationary.  Not surprisingly, this turns out to be true of H2 as well.  In
particular 

\Prop \newcount\Sphere\Sphere=\the\ThisProp If $\ladder\varphi\lambda\Lambda$ 
is an H2-sequence, then for any stationary set $S$, $\{\lambda\in S\,:\, S\cap\fsubl\neq\eset\}$
is stationary.\endProp

There is a similar proposition for H1. 
What is more surprising than this proposition is the following:

\Prop If
$\ladder\varphi\lambda\Lambda$ is an H2-sequence, then for any set $C\subseteq\omega_1$ which
contains a club, the set $\{\lambda\in C\,:\,C\cap\fsubl\neq\eset\}$ contains a club.
\showProp Let $\ladder\varphi\lambda\Lambda$ be an H2-sequence and suppose $C\subseteq
\omega_1$ is a set containing a club for which the set in question does not contain a club.
Then the complement of this set, namely $\{\lambda\,:\, \lambda\notin C \,\vee\,
C\cap\fsubl=\eset\}$, is stationary.  We may intersect this set with $C$ to get a stationary
set:  $A\Odef\{\lambda\in C\, :\,C\cap\fsubl=\eset\}$.  Since the $\fsubl$ form an H2-sequence,
there is a $\lambda\in A$ such that $A\cap\fsubl\neq\eset$.  But $A\subseteq C$ implies
$\fsubl\cap C\neq\eset$ which by the definition of $A$ implies $\lambda\notin A$, a
contradiction.  \endProp

 As an aside, Proposition \the\Sphere\  indicates there is an intermediary
between H2 and weak-H2, namely:

\DDef{} Not-as-weak-H2 is the statement that there is a ladder system
$\ladder\varphi\lambda\Lambda$ such that for all functions $f\,:\,\omega_1\to\omega$, the set
$\{\lambda\,:\,f^{-1}(\lambda) \cap\fsubl\neq\eset\}$ is stationary.  \endDDef

Again, there is a similar statement for H1.


\Section{The $\aleph_1$ Gap Theorem}

In this section, we prove a strengthening of the $\aleph_0$ Gap Theorem.  This theorem is not 
stated as a consistency result.  It is a construction that occurs in ZFC.  However, to show that 
the constructed object is in fact a large collection of gaps, and in particular to prove the
corollary corresponding to that following the $\aleph_0$ Gap Theorem, we use hypothesis H2.  
Towards this corollary, we first prove the following lemma which shows that under H2 a weak 
condition on a pregap makes it a gap.  

\Lemma(H2 and Gaps Lemma.)
\LocalLLemma=\the\ThisProp
  Assume H2 holds for the ladder system $\omladder$.  Let $\langle
A,\,B\rangle$ be a pregap satisfying $A_{\alpha}\cap B_{\alpha}=\eset$ for all $\alpha$ and
$$(\exists\eta )(\exists\text{ stationary }S\subseteq\Lambda )(\forall \lambda\in
S)(\forall\beta\in\varphi_{\lambda}\Sminus(\eta +1))\,\, A_{\lambda}\cap B_{\beta}\neq\emptyset
.$$ Then $\langle A,\,B\rangle$ is a gap.

The appearance of $\eta$ in this lemma is for technical reasons that will be clear in its
application after the proof of the $\aleph_1$ Gap Theorem.  The proof of the lemma is better
understood ignoring $\eta$.

\showLemma Suppose $D\subseteq\omega$ separated $\langle A,B\rangle$.  Then there is an
$n\in\omega$ and a stationary set $T\subseteq S\Sminus(\eta +1)$ such that $(\forall\lambda
,\lambda'\in T)$ $$\matrix A_{\lambda}\Sminus n\subseteq D&(B_{\lambda}\Sminus n)\cap
D=\emptyset\cr A_{\lambda}\cap n=A_{\lambda'}\cap n&B_{\lambda}\cap n=B_{\lambda'}\cap
n.\endmatrix $$ By H2, $(\exists\lambda\in T)\,\varphi_{\lambda}\cap T\neq\emptyset$.  Let
$\beta\in\varphi_{\lambda}\cap T$.  Since $\beta\in\varphi_{\lambda}$, $A_{\lambda}\cap
B_{\beta}\neq\emptyset$, while $\beta$ and $\lambda$ being in $T$ implies $A_{\lambda}\cap
B_{\beta}=\emptyset$, by the previous displayed equations.  This is a contradiction, and hence
no such $D$ exists.  \endLemma

We now state and prove the title theorem of this section.  
\Thm {\bf The} $\aleph_1$ {\bf Gap Theorem.\quad }Let the following objects be given.  

\leftskip.1truein 
\item{(G1)} A tower $T=\langle T_{\alpha}:\alpha\in\omega_1\rangle$ with $T_0=\eset$.
\item{(G2)} A ladder system $\ladder{\varphi}\lambda\Lambda$.
\item{(G3)} A collection of disjoint stationary subsets of $\Lambda$:  $\langle\Sxh:\eta
<\xi
<\omega_1\rangle$.  We shall also assume that each $\lambda\in\Lambda$ is in some $\Sxh$ and 
that \hfil\newline $(\lambda\in\Sxh)\implies (\eta <\xi <\lambda )$.  

\leftskip0truein Then there is a collection $\langle\prxa A:\xi <\alpha <\omega_1\rangle$ of
subsets of $\omega$ satisfying:

\leftskip.2truein 
\item{($\aleph_1$ 1)} $(\forall\xi\in\omega_1)\,\langle\prxa A:\xi <\alpha <\omega_1\rangle$ is 
a  subtower of $T$.  
\item{($\aleph_1$ 2)} $(\forall\alpha\in\omega_1)\,T_{\alpha}=\bigcup \{\prxa A:\xi <\alpha 
\}$, a disjoint union.  
\item{($\aleph_1$ 3)} $(\forall\eta <\xi <\omega_1)(\forall\lambda\in\Sxh 
)(\forall\beta\in\fsubl.\beta >\eta )A(\xi, \lambda) \cap A(\eta, \beta)\neq\eset$.  
\item{($\aleph_1$ 4)} $(\forall\beta <\alpha <\omega_1)$ the set $\{\xi :\prxb A\nsubseteq\prxa 
A\}$ is finite.  

\leftskip0truein

Thinking of the indexing of $\prxa A$ as $A(\text{column, row})$, we can visualize the result as
an $\omega_1\times\omega_1$ upper triangular matrix with $\pr00$ at the lower left corner.  This
follows the intention of building $\omega_1$-towers (the columns) while keeping each row a
countable disjoint collection of sets whose union is $T_{\alpha}$.

These conditions on the matrix will be satisfied by a recursive construction of sets at level,
$\gamma$ ($\lambda$ when working with limit stages).  Consequently we will frequently refer to
the ``restrictions'' of ($\aleph_1$~1--4) as the {\it induction hypotheses}, which are obtained
by replacing $\omega_1$ by $\gamma$ (or $\lambda$) and quantifying over the sets constructed to
that point in the proof.  When no confusion can result, we shall refer to these restrictions as
($\aleph_1$~1), etc.

There is a corollary to this theorem corresponding to the one after the $\aleph_0$ Gap Theorem.
However, the corollary is stated in terms of an independence result and its proof is more 
involved, so we delay its consideration until after the construction.  

\showThm 

Assume the objects in (G1)--(G3) have been fixed.  We construct sets $A(\xi,\,\alpha)$ with the
convention that $\prxa A=\eset$ for $\xi\geq\alpha$.  We first establish the following lemma 
which provides a convenient equivalent formulation of ($\aleph_1~4$).  

\Lemma Suppose that sets $\prxa A$ for $\xi <\alpha <\gamma$ satisfy the restrictions of 
$(\aleph_1~1)$ and $(\aleph_1~2)$.  Fix $\beta <\alpha <\gamma$.  Then the following are 
equivalent:  \LocalLemma=\the\ThisProp 
\item{(A)} The set $U_1\Odef\{\xi :\prxb A\nsubseteq\prxa A\}$ is finite.  
\item{(B)} The set $U_2\Odef\{\xi :\prxa A\cap T_{\beta}\neq\prxb A\}$ is finite.  

\showLemma 

\noindent(A) $\implies$ (B):  Assume $U_1$ is finite.  If $\xi\in U_2$ then at least one of the
following three conditions holds of $\xi$:  
\item{(i)} $\prxb A\nsubseteq T_{\alpha}$, or 
\item{(ii)} The set $(\prxb A\Sminus\prxa A)\cap T_{\alpha}$ is non-empty, or 
\item{(iii)} The set $\prxa A\cap T_{\beta}\Sminus\prxb A$ is non-empty.  

It is sufficient to show that there are only finitely many $\xi$ satisfying each of (i), (ii)
and (iii).  This is true for (i) by our assumption that $(\aleph_1~1)$ and $(\aleph_1~ 2)$ hold.
Next, if $\xi$ satisfies (ii) then $\prxb A\nsubseteq\prxa A$, and so $ \xi\in U_1$ which is
assumed to be finite.

If $\xi$ satisfies (iii) then by $(\aleph_1~2)$ for level $\beta$, $ (\exists\eta )\prxa
A\cap\prhb A\neq\emptyset$.  Let $$U_2^{+}(\eta )\Odef\{\xi\,:\,\prxa A\cap\prhb
A\neq\emptyset\land \eta\neq\xi \}.$$ Then the set of $\xi$ satisfying (iii) is equal to
$\bigcup_{\eta <\gamma}U_2^{+}(\eta )$.  By $(\aleph_1~2)$ for level $\alpha$, if $U_2^{+}(\eta
)\neq\emptyset$ then $\eta\in U_1$.  Since we are assuming that $U_1$ is finite, it must be the
case that for only finitely many $\eta$ is the set $U_2^{+}(\eta )$ non-empty.  It remains to
show that each $U_2^{+}(\eta )$ is finite.  By $(\aleph_ 1~2)$, the $\prxa A$ are disjoint.  So,
if for some $\eta$ the set $U_2^{+}(\eta )$ were infinite, then $\prhb A\nsubseteq^{*}A(\eta
,\,\alpha )$, contradicting the restriction of $ (\aleph_1~1)$ for column $\eta$.  Hence only
finitely many $\xi$ satisfy (iii).

\noindent(B) $\implies$ (A):  This is immediate since $\prxb A\nsubseteq\prxa A$ implies that
$\prxa A\cap T_{\beta}\neq\prxb A$.  \endLemma

Consequently, we may assume (B) holds, but need only prove that (A) is maintained.

\noindent {\it The Construction of the Matrix.  }

The stage $\gamma =0$ is trivial.

\noindent {\it Successor Stage:} $\,\gamma =\alpha +1$.

Assume that for $\xi <\beta <\gamma$, sets $\prxb A$ have been constructed satisfying the
induction hypotheses.  Let $\langle B(\xi ):\xi <\gamma\rangle$ be a partition of
$T_{\gamma}\Sminus T_{\alpha}$ into infinite disjoint sets.  Define $$A(\xi ,\,\gamma
)\Odef(\prxa A\cap T_{\gamma})\cup B(\xi )\text{ for each $ \xi <\gamma$.}$$ (Recall that by
convention $A\pr\alpha\alpha =\eset$.)

Clearly ($\aleph_1$ 1) and ($\aleph_1$ 2) now hold.  The condition in ($\aleph_1$ 3) will be
satisfied at limit stages.  So it remains to check ($\aleph_1$ 4), in particular that
$$(\forall\beta <\gamma )\,\,\{\xi :\prxb A\nsubseteq A(\xi ,\gamma ) \}\text{ is finite.}$$

{\it Case 1:} $\beta =\alpha${\it .\/} Since $T_{\alpha}\Sminus T_{\gamma}$ is finite and since
the rows are disjoint families, it follows from the definition of $A(\xi ,\gamma )$ that the
desired set is finite.  In particular, it is contained in the set $\{\xi :\prxa A\cap
(T_{\alpha}\Sminus T_{\gamma})\neq\emptyset \}$.

{\it Case 2:} $\beta <\alpha${\it .\/} Taking the contrapositive of $$\left[\,\prxb
A\subseteq\prxa A\,\land\,\prxa A\subseteq A(\xi ,\gamma )\,\right]\implies\prxb A\subseteq
A(\xi ,\gamma )$$ we have $$\{\xi :\prxb A\nsubseteq A(\xi ,\gamma )\}\subseteq \{\xi :\prxb
A\nsubseteq\prxa A\}\cup \{\xi :\prxa A\nsubseteq A(\xi ,\gamma ) \}.$$ Both of the sets on the
right are finite, the first by induction hypothesis ($\aleph_1~4$) and the second by Case 1.

This establishes the preservation of ($\aleph_1$ 4) through the successor stage and so completes
this part of the construction.  We now address the arduous

\noindent {\it Limit Stage:} $\lambda$.

Assume that $\lambda\in S(\xi_0,\,\eta_0)$ for some $\eta_0<\xi_ 0<\lambda$, and that $\prxa A$
are constructed for $\xi <\alpha <\lambda$ satisfying the induction hypotheses.  Fix bijections,
$f:\omega\to\lambda$, and $g:\omega\to (\lambda\times 2)\cup T_{\lambda}$.

The following is a brief description of the construction which is to follow.  The method is
expressed in the notation and ideas of forcing.  However, all objects involved are countable so
no new generic objects are needed to obtain the result.  Nonetheless, future extensions of this
theorem may take advantage of this methodology.

We describe six different properties obtainable by objects of the form $\ot rsn$ where
$n\in\omega$, and $r$ and $s$ are finite partial functions on $\lambda$ with codomains $\omega$
and $\twopre n$ respectively.  Those tuples satisfying these properties will be called {\it
conditions}.

The function $s(\xi )$ is a finite approximation of (the characteristic function of) the set
$A\pr\xi\lambda$.  If $\xi\in dom(s)$ and $\beta \in dom(r)$, we consider this a promise to
satisfy $A(\xi ,\,\beta )\Sminus n\subseteq A(\xi ,\,\lambda )$.  The value of $r(\beta
)=m\in\omega$ will be a promise to satisfy ($\aleph_1$ 4) between rows $\lambda$ and $\beta$
``above'' $m$ (really, to satisfy the condition for $\xi$ whose $f$-preimage is greater than
$m$).  Thus if $f^{- 1}(\xi )>r(\beta )$ then we promise to have $A(\xi ,\,\beta )\subseteq
A(\xi ,\,\lambda )$.  $n$ is redundant, but convenient to have explicit as it is frequently
referenced.

The idea of forcing is implicit in this description.  Here is an informal list of what a
condition forces.

\leftskip 0.1truein 
\item {(F1)} A number, $i\in T_{\lambda}$, is explicitly forced into some $A(\xi ,\,\lambda )$ 
by $\ot rsn$ if $\xi\in dom(s)$ and $s(\xi )(i)=1$.  

\item {(F2)} A number, $i\in T_\lambda$ may be implicitly forced into some $A(\xi ,\,\lambda )$
by $\ot rsn$ due to the ``almost containment'' for the tower:  $i\geq n$, $\xi\in dom(s)$ and
there is a $\beta\in dom(r)$ with $ i\in A(\xi ,\,\beta )$, or

\item {(F3)} A number $i\in T_\lambda$ can also be implicitly forced into $A(\xi ,\,\lambda )$
due to $(\aleph_1~4)$:  there is a $\beta\in dom(r)$ with $f^{-1}(\xi )>r(\beta )$ and $i\in
A(\xi ,\,\beta )$.

\item {(F4)} A number $i\in T_{\lambda}$ is forced out of an $A(\xi ,\,\lambda )$ by $\ot rsn$
explicitly if $s(\xi )(i)=0$ or implicitly just in case it is forced into $A(\eta ,\,\lambda )$
for some $\eta \neq\xi$.  Of course, $i$ is forced out of all sets if $i\notin T_{\lambda}$.

\item{(F5)} Hence, $\ot rsn$ forces $\prxb A\Sminus n\subseteq A(\xi ,\,\lambda )$ whenever
$\xi\in dom(s)$ and $\beta\in dom(r)$.

\item{(F6)} To ensure $(\aleph_1~3)$ is satisfied, $\ot rsn$ will explicitly force that
$A\pr{\eta_0}\delta\cap A\pr{\xi_0}\lambda$ is non-empty for all $\delta$ satisfying
$\delta\in\fsubl\cap\sup^{+}dom(r)\land\delta >\eta_0$.

\leftskip 0 truein \noindent Most of the definition of \bigP\ below can be seen as consistency
requirements for this ``forcing''.

A partial order called {\it extension\/} will be defined on conditions with the intuition being
that an extension contains more information about the sets being constructed.  A four part
extension lemma is proved with the following implications:  (E1) permits $n$ to be incremented
and is essentially a service lemma for the following parts.  (E2) and (E3) permit the extension
of the domains of $s$ and $r$, respectively, by an element.  (E4) permits the addition of an
arbitrary element of $T_{\lambda}$ into some $A\pr\xi\lambda$.

A recursive definition is given starting with the initial condition $\langle\emptyset
,\,\{\langle\xi_0,\,\emptyset\}\rangle ,\,0\rangle$, to which is applied the appropriate extension
lemma which is dictated by the type of $g(m)$ where $m$ is the stage of the definition.

This generates a chain of objects from which the $A(\xi ,\lambda )$ are derived.  $g$ is a
bookkeeping function that ensures all the desired properties are obtained.  This ends the
description of the proof mentioned above.

For a function $s(\xi )$ mapping into 2, recall the notation that $$\sbar\xi\Odef(s(\xi
))^{-1}\{1\}.$$

\DDef{\noexpand\bigP}Define the set $$\align\text{\bigP}\subseteq \{\ot
rsn:\,\,&dom(r),\,dom(s)\in [ \lambda ]^{<\aleph_0}\,\land\,\ninom\,\land\\
&ran(r)\subseteq\omega\,\land\,ran(s)\subseteq{}^n2\}\endalign$$ as follows.

$\ot rsn\in$ \bigP\ if and only if the following requirements are satisfied.  Note that the
parenthetical statements are meant as explanation, not as part of the definition.

\leftskip.1truein
\item{(\bigP1)} $\langle\sbar\xi :\xi\in dom(s)\rangle$ is a collection of
pairwise disjoint subsets of $T_{\lambda}$.  ($s$ builds the final sets.  See 
(F1).)
\item{(\bigP2)} $(\forall\alpha\in dom(r))(\forall\xi\in dom(s))(\forall m>r(
\alpha ))$
\itemitem{(\bigP2A)} $\xi =f(m)\implies\prxa A\cap n\subseteq\sbar
\xi$,
(($\aleph_1$ 4) satisfied with row $\alpha$.  See (F3).)
\itemitem{(\bigP2B)} $\xi\neq f(m)\implies A\pr{f(m)}\alpha\cap\sbar
\xi =\eset$.
(Allows $f(m)$ to be added to the domain of $s$.)
\itemitem{(\bigP2C)} $A\pr{f(m)}\alpha\subseteq T_{\lambda}$.
\item{(\bigP3)} $(\forall\alpha ,\,\beta\in dom(r))(\forall\xi\in dom(s))
(\forall m>r(\alpha ))$ \hfill \newline 
$\xi\neq f(m)\implies\prxb A\cap A\pr{f(m)}\alpha\Sminus n=\eset$.  (See (F5).)
\item{(\bigP4)} $(\forall\alpha ,\,\beta\in dom(r))(\forall m_0>r(\alpha 
))(\forall m_1>r(\beta ))$ \hfill \newline 
 $m_0\neq m_1\implies A\pr{f(m_0)}\alpha\cap A\pr{f(m_1)}\beta =\eset$.
\hfill \newline 
(When $f(m_0)$ and $f(m_1)$ are added to the domain of $s$, we need to 
meet ($\aleph_1$ ~4) at $\alpha$, $\beta$ respectively while maintaining
disjointness of $\sbar\xi$.)
\item{(\bigP5)} $(\forall\xi\neq\xi'\in dom(s))(\forall\beta ,\beta'\in d
om(r))\,\,\prxb A\cap A(\xi',\,\beta')\Sminus n=\eset$. (Needed to allow almost
containment to be satisfied.  See (F5).)
\item{(\bigP6)} $\xi_0\in dom(s)$ and $(\forall\delta\in\fsubl\cap\sup^{
+}dom(r).\delta >\eta_0)\,A\pr{\eta_0}\delta\cap\sbar{\xi_0}\neq\emptyset 
.$  (This will ensure
condition ($\aleph_1$ 3) is satisfied.  See (F6).)
\endDDef

\leftskip0truein
\DDef{$\leq$ on \noexpand\bigP} If $p=\ot rsn$ and $p'=\otp rsn$ are elements of
\bigP,
say $p'$ {\it extends} $p$ and write $p'\leq p$ if and only if the following 
five conditions are satisfied:
\item{(e1)} $r'\supseteq r$;
\item{(e2)} $n'\geq n$;
\item{(e3)} $dom(s')\supseteq dom(s)$;
\item{(e4)} $(\forall\xi\in dom(s))\,s'(\xi )\supseteq s(\xi )$;
\item{(e5)} $(\forall\xi\in dom(s))(\forall\beta\in dom(r))\,\sbarp\xi\supseteq 
(\prxb A\cap (n'\Sminus n))$.  (See (F5).)
\endDDef

We now state the Extension Lemmas.  For all these statements, let
$p=\ot rsn\in$ \bigP.
\Lemma(E1).  (The Simple Extension Lemma.)  There is an $s'$ such that 
\item{(a)} $\langle r,\,s',\,n+1\rangle\in$ \bigP,
\item{(b)} $\langle r,\,s',\,n+1\rangle\leq p$,
\item{(c)} $dom(s')=dom(s)$.
\Lemma(E2).  (Extending the domain of $s$.)  For any $\xi\in\lambda$ 
there are an $s'$ and an $n'$ such that
\item{(a)} $\langle r,\,s',\,n'\rangle\in$ \bigP,
\item{(b)} $\langle r,\,s',\,n'\rangle\leq p$,
\item{(c)} $dom(s')=dom(s)\cup \{\xi \}$.
\Lemma(E3).  (Extending the domain of $r$.)  For any $\beta\in\lambda$ 
there are an $s'$, an $n'$ and an $m'$ such that
\item{(a)} $\langle r\cup \{\langle\beta ,\,m'\rangle \},\,s',\,n'
\rangle\in$ \bigP,
\item{(b)} $\langle r\cup \{\langle\beta ,\,m'\rangle \},\,s',\,n'
\rangle\leq p$,
\item{(c)} $dom(s')=dom(s)$.
\Lemma(E4).  (Adding an element of $T_{\lambda}$.)  For any $i\in 
T_{\lambda}$ 
there are an $s'$ and an $n'$ such that
\item{(a)} $\langle r,\,s',\,n'\rangle\in$ \bigP,
\item{(b)} $\langle r,\,s',\,n'\rangle\leq p$,
\item{(c)} $(\exists\xi\in dom(s'))\,i\in\sbarp\xi$,
\item{(d)} $\vert dom(s')\Sminus dom(s)\vert\leq 1$.

\showLemma(E1).

Let $\ot rsn\in$ \bigP\ be given.  To show the claim, we must extend each function $s(\xi )$
by one place such that the resulting triple $\langle r,\,s',\,n+1\rangle$ satisfies (\bigP1)
through (\bigP6).  Notice that since $r$ and $dom(s)$ do not change, (\bigP2C) and (\bigP3)
through (\bigP6) will necessarily be satisfied.  By (\bigP1), $n$ should belong to at most one
$\sbar\xi$.  We perform a minimal extension to satisfy (\bigP2A) and the definition of
extension, (e5).  That is, define $s'(\xi ) (n)\Odef0$ unless there is an $\alpha\in dom(r)$
with $n\in\prxa A$ in which case $s'(\xi )(n)\Odef1$.  If this new triple is in \bigP, it
immediately satisfies the conditions to extend $\ot rsn$.  So, we show

\Claim$\langle r,\,s',\,n+1\rangle\in$ \bigP.  
\showClaim We first show that only one such $\xi$ can satisfy this second requirement.  That is 
if $\ot rsn\in$ \bigP, then for at most one $\xi$ is there an $\alpha\in dom(r)$ with $n 
\in\prxa A$.  But this is exactly what (\bigP5) states for $\ot rsn$.  Hence $\langle 
r,\,s',\,n+1\rangle$ satisfies (\bigP1).  

The condition defining $\langle r,\,s',\,n+1\rangle$ gives (\bigP2A) immediately.  If
(\bigP2B) failed, there would be an $\alpha\in dom(r)$, a $\xi\in dom(s)$, and an $m>r(\alpha
)$ such that $\xi\neq f(m)$ and $n\in\prxa A\cap A(f(m),\,\alpha )$.  But this contradicts the
fact that (\bigP3) held for $\ot rsn$.  This completes the proof that $\langle
r,\,s',\,n+1\rangle\in$ \bigP.  \endClaim

This completes the proof of (E1).

Before continuing, we state and prove a lemma necessary for the remaining proofs.  While (E1)
stated that $n$ could be increased, the following lemma shows there is an $n'$ to which $n$
can be increased to meet the other conditions in the definition of condition.

\Lemma\ \LocalLLLemma=\the\ThisProp\ Given $\langle r,s,n\rangle\in$ \bigP, and $\xi\notin
dom(s)$ there is an $n'>n$ satisfying the following two conditions:
\item{(a)} $(\forall\alpha\neq\beta\in dom(r))(\forall m>r(\alpha ))\,\prxa A\cap
A(f(m),\,\beta )\Sminus n'=\emptyset .$
\item{(b)} $(\forall\eta\in dom(s))(\forall\alpha\neq\beta\in dom(r))\,\prxa A\cap\prhb
A\Sminus n'=\emptyset .$

\showLemma Since (a) and (b) are preserved as $n'$ grows and since there are only finitely
many triples $\eta\in dom(s),\,\alpha\neq \beta\in dom(r)$, it suffices to show such an $n'$
exists for an arbitrary such triple.  For (b), notice that as $\xi\notin dom(s)$ we have
$\xi\neq\eta$.  Since $\prxa A$ and $\prhb A$ are almost disjoint, (b) follows.

For (a) there are two cases depending on the order of $\alpha$ and $ \beta$.

\noindent$\beta >\alpha$:  By induction hypothesis and Lemma \the\LocalLemma,
$\{\mu :A\pr\mu\beta\cap T_{\alpha}\neq A\pr\mu\alpha \}$ is finite.  So there is an $m_0$
such that for all $m$ $$\align m>m_0\land\xi\neq f(m)&\implies A\pr{f(m)}\beta\cap T_{
\alpha}=A\pr{f(m)}\alpha\\ &\implies A\pr{f(m)}\beta\cap\prxa A=\eset.\endalign$$ But there
are only finitely many $m$ with $r(\beta )<m\leq m_0$ while $\xi\neq f(m)$ implies
$A\pr{f(m)}\beta$ and $\prxa A$ are almost disjoint.

\noindent$\beta <\alpha$:  Similarly to the previous case, the set $\{\mu
:A\pr\mu\beta\nsubseteq A\pr\mu\alpha \}$ is finite.  So there is an $m_0$ such that for all
$m$ $$\align m>m_0\land\xi\neq f(m)&\implies A\pr{f(m)}\beta\subseteq A\pr{f(m)}\alpha\\
&\implies A\pr{f(m)}\beta\cap\prxa A=\eset.\endalign$$ The case is completed as above.
\endLemma

\showLemma(E2).

Fix $\ot rsn\in$ \bigP\ and $\xi\in\lambda\Sminus dom(s)$.  We wish to add $\xi$ to the domain
of $s$.  This is done in two steps, first by extending the functions $s(\eta )$ for the
$\eta\in dom(s)$ (i.e., increasing $n$) and then by defining $s'(\xi )$.

\noindent {\it Step 1:  \/} Fix $n'$ as in Lemma \the\LocalLLLemma\ for $\ot
rsn$.  By iterated application of (E1), we may extend $\ot rsn$ to $\langle
r,\,s^{\prime\prime},\, n'\rangle\in$ \bigP.  This ensures that (\bigP3) and (\bigP5) will be
satisfied by the new condition.

\noindent {\it Step 2:\/} We now add $\xi$ to $dom(s)$ conforming to condition (\bigP2A).  For
$\eta\in dom(s)$, let $s'(\eta )\Odef s^{\prime\prime}(\eta )$.  Define $s'(\xi )$ by
$$\sbarp\xi\Odef\{i<n'\,:\,(\ex\alpha\in dom(r))\,\,i\in\prxa A\, \land\,f^{-1}(\xi )>r(\alpha
)\}.$$ 

\Claim$\langle r,\,s',\,n'\rangle\in$ \bigP.  
\showClaim Sketch.

The construction assures that (\bigP2A), (\bigP3) and (\bigP5) are satisfied.  (\bigP2C),
(\bigP4) and (\bigP6) persist from $\langle r,\,s^{\prime\prime},\,n'\rangle$.  If (\bigP1)
fails for $\langle r,\,s',\,n'\rangle$, then (\bigP2B) would not hold for $\langle
r,\,s^{\prime\prime},\,n'\rangle$.  If (\bigP2B) fails for $\langle r,\,s',\,n'\rangle$, then
(\bigP4) would not hold for $\langle r,\,s^{\prime\prime},\,n'\rangle$.  \endClaim

Hence $\langle r,\,s',\,n'\rangle\in$ \bigP.  The other desired conditions for (E2) follow
from the construction.  This completes the proof of (E2).  \endLemma 

\showLemma(E3).

\noindent {\it Part 1:\/} Let $\beta\in\lambda\Sminus dom(r)$.  We wish to add $\beta$ to the
domain of $r$.  The initial part of the construction is more complicated if $\beta
>\max(\sup(dom(r)),\eta_0)$.  If this is {\it not\/} the case, we may let $\langle
r,\,s^{\prime\prime},\,n^{\prime\prime}\rangle\Odef\ot r sn$ and skip to {\it Part 2\/} of the
construction.

But suppose $\beta >\max(\sup(dom(r)),\eta_0)$.  Recall that $\lambda\in\Sxhz$ and by
(\bigP6), $\xi_0\in dom(s)$.  We must ensure that (\bigP6) is satisfied by the extension, in
particular for each $\delta\in (\fsubl\cap (\beta +1)\Sminus(\eta_0+1))$, $\,A\pr{\eta_
0}\delta\cap\sbarp{\xi_0}\neq\eset$.

Consider an arbitrary $\delta\in (\fsubl\cap (\beta +1))\Sminus\max
(\sup^{+}(dom(r),\eta_0+1))$, a finite set.  We wish to find an ``unrestricted'' number in
$A(\eta_0,\,\delta )$ to force into $ \bar {s}'(\xi_0)$.  Say that $i\in
A\pr{\eta_0}\delta\Sminus n$ is {\it restricted\/} (i.e., is already forced into another set)
if either of the following conditions are met:

\item{(1)} $(\ex\alpha\in dom(r))(\ex m>r(\alpha ))\,\,i\in A(f(m ),\,\alpha )$, or
\item{(2)} $(\ex\alpha\in dom(r))(\ex\xi\in dom(s))\,\,i\in\prxa A$.

Since $\delta >\sup^{+}(dom(r))$, the set $A\pr{\eta_0}\delta\Sminus \bigcup
\{T_{\alpha}:\alpha\in dom(r)\}$ is infinite.  (Actually, we also use the following facts:
($\aleph_1$ 4) is satisfied as an induction hypothesis, sets in different columns are almost
disjoint and $A\pr{\eta_0}\delta\Sminus A\pr{\eta_0}\alpha$ is infinite for each $\alpha\in
dom(r)$.)  This means there are arbitrarily large $ i$ which are not restricted.

Perform the following extension process for each relevant $\delta$ in turn.  Let $i_{\delta}$
be the least $i\geq n$ such that $$i_{\delta}\in\left(A\pr{\eta_0}\delta\Sminus\bigcup
\{T_{\alpha}
:\,\alpha\in dom(r)\}\right)\cap T_{\lambda}.$$ (That is, $i_{\delta}$ is not restricted.)
Repeatedly apply (E1) to get $\langle r,\,s^{\prime\prime},\,i_{\delta}\rangle\leq\ot rsn$
with $ dom(s^{\prime\prime})=dom(s)$.  Extend this in the obvious minimal way to satisfy
$i_{\delta}\in\bar {s}^{\prime \prime}(\xi_0)$.  Notice that as $i_{\delta}$ is not
restricted, this is a condition extending $\ot rsn$.  Let the result of this finite iteration
be called $\langle r,\,s^{\prime\prime},\,n^{\prime\prime}\rangle$.

\noindent {\it Part 2:\/} Assume we have a condition $\langle r,\,
s^{\prime\prime},\,n^{\prime\prime}\rangle$ satisfying (\bigP6) for $\delta\in\fsubl\cap
(\beta +1)\Sminus(\eta_ 0+1)$.  We now extend $s^{\prime\prime}$ and $n^{\prime\prime}$ to
satisfy (\bigP3) and (\bigP5).  By a slight variation of the Lemma
\the\LocalLLLemma, there is an $n'\geq n^{\prime\prime}$ such that
\item{(a)} If $\alpha\in dom(r)$, $m>r(\alpha )$ and $\xi\in dom( s)$ satisfy $\xi\neq f(m)$
then $\prxb A\cap A\pr{f(m)}\alpha\Sminus n'=\eset$.
\item{(b)} If $\xi\neq\xi'\in dom(s^{\prime\prime})$ and $\beta'\in dom(r)$ then $\prxb A\cap
A\pr{\xi'}\beta'\Sminus n'=\eset$.  Iteratively applying (E1), obtain $\langle
r,\,s',\,n'\rangle\leq \langle r,\,s^{\prime\prime},\,n^{\prime\prime}\rangle$.

\noindent {\it Part 3:\/} Defining $r(\beta )$.

\Claim There is an $m'\in\omega$ satisfying the following three conditions:
\item{(1)} $(\forall m>m')\,A\pr{f(m)}\beta\cap n'=\eset$.
\item{(2)} $(\forall\alpha\in dom(r))(\forall m>m')(\forall l>r(\alpha ))$ \hfill \newline
$l\neq m\implies A\pr{f(l)}\alpha\cap A\pr{f(m)}\beta =\eset$.
\item{(3)} $(\forall m>m')\,A\pr{f(m)}\beta\subseteq T_{\lambda}$.  

\showClaim Since there are
only finitely many requirements on $m'$, and since the properties of interest are preserved
upward (for larger $m'$), we may consider each one separately.  (1) follows from the fact that
$\{\prhb A:\,\eta <\beta \}$ is a disjoint family while (3) follows because $T_{\beta}\Sminus
T_{\lambda}$ is finite.

For (2) fix $\alpha\in dom(r)$.  If $\alpha <\beta$, take $m'>\max(f^{-1}\{\xi :\,\prxb A\cap
T_{\alpha}\neq\prxa A\})$, a finite set by hypothesis.  If $\alpha >\beta$, let
$m'>\max(f^{-1}\{\xi :\,\prxb A\nsubseteq\prxa A\}$, again finite by hypothesis.  In either
case, if $\xi\neq f(m)$ and $m>m'$ then $A\pr{f(m)}\beta\cap\prxa A=\eset$.  \endClaim

Fix $m'$ as in the previous claim.  

\Claim$\langle r\cup \{(\beta,\,m' )\},\,s',\,n'\rangle\in$ \bigP.  
\showClaim All the work has 
been done in the previous lemmas and claims.  For example, (\bigP1) holds for $\langle 
r,\,s',\,n'\rangle$ and hence holds for $\langle r\cup \{(m',\,\beta )\},\,s',\,n'\rangle$.  
(\bigP2) holds by choice of $m'$, condition (3), as does (\bigP3) when $\beta$ is $\alpha$.  
When $\beta$ instantiates $\beta$ in (\bigP3), the choice of $n'$, part (2), ensures the 
condition is met.  (\bigP4) holds by choice of $m'$, part (3), (\bigP5) by choice of $n'$ part 
(2), and finally (\bigP6) holds because of the construction of $\langle 
r,\,s^{\prime\prime},\,n^{\prime\prime}\rangle$.  \endClaim 

The other properties of this condition are clear from the construction and so this completes
the proof of (E3).

\showLemma(E4).  Fix $i\in T_{\lambda}$.  We wish to find an $s'$ and an $n'$ such that there
is a $\xi\in dom(s')$ with $i\in\sbarp\xi$ and with $dom(s')$ having at most 1 new element.
The conclusion is immediate if $i$ is in $\sbar\xi$ for some $\xi$ in $dom(s)$, so assume
otherwise.  We examine three cases.

{\it Case 1:} $\ot rsn$ implicitly ``forces'' $i\in A\pr\xi\lambda$ for some $\xi$.  See (F2)
and (F3).  Either:
\item{(1)} There is an $\alpha\in dom(r)$ and an $m>r(\alpha )$ such that $i\in
A(f(m),\,\alpha )$, or
\item{(2)} $i>n$ and there is an $\alpha\in dom(r)$ and $\xi\in dom(s)$ such that $i\in\prxa
A$.

Either way, apply (E1) and (E2) to ensure that $\langle r,\,s',\, n'\rangle\leq\ot rsn$,
$n'>i$ and $\xi\in dom(s')$.  Then the conclusion holds by (\bigP2A) applied to $\langle
r,\,s',\,n'\rangle$ or by (e5) according to whether (1) or (2) holds, respectively.

{\it Case 2:\/} Case 1 fails and $i<n$.  Then we may fix $\xi\notin dom(s)$ and define a new
condition that satisfies $i\in\sbarp\xi$.  This case uses (\bigP4) on $\ot rsn$ to ensure the
result holds.

{\it Case 3:\/} Case 1 fails and $i\geq n$.  Extend by (E1) to $\langle r,\,s',\,i\rangle$ and
put $i$ in any $\sbarp{\xi}$.

These exhaust the cases need to establish (E4) and complete the proof of the Extension Lemmas.
\endLemma

\noindent {\it Completion of the Construction.}

We recursively define a sequence of conditions $p_m$ for $m\in\omega$.  Let
$p_0\Odef\langle\emptyset ,\,\{\langle\xi_0,\,0\rangle \},\,0 \rangle$ and note that $p_0\in$
\bigP.  Recall that $g:\omega\to (\lambda\times 2)\cup T_{\lambda}$ is a bijection.  Suppose
that $p_m=\otss mrsn$ is defined.  We define $p_{m+1}\leq p_m$ according to the value of
$g(m)$ as follows:

{\it Case 1:} $g(m)=\op\xi 0$ for some $\xi <\lambda$.  Put $\xi$ into the domain of $s$.
Apply Extension Lemma (E2) to $p_m$ and $\xi$ to get $p_{m+1}=\otss{m+1}rsn\leq p_m$ with
$r_{m+1}=r_m$ and $\xi\in dom(s_{m+1})$.

{\it Case 2:} $g(m)=\op\beta 1$ for some $\beta <\lambda$.  Put $ \beta$ into the domain of
$r$.  Apply Extension Lemma (E3) to get $p_{m+1}=\otss{m+1}rsn\leq p_m$ with $\beta\in
dom(r_{m+1})$.

{\it Case 3:} $g(m)\in T_{\lambda}$.  Ensure $g(m)$ is in some $\sbar \xi$.  Apply Extension
Lemma (E4) to obtain $p_{m+1}=\otss{m+1}rsn\leq p_ m$ with $g(m)\in\bar {s}_{m+1}(\xi )$ for
some $\xi\in dom(s_{m+1})$.

This completes the definition of $p_m$ and we are ready to define the $A\pr\xi\lambda$.  Fix
$\xi <\lambda$.  Since $g$ is a bijection, there is an $l\in\omega$ with $g(l)=\op\xi 0$.  By
the definition of $p_m$, if $m>l$ then $\xi\in dom(s_m)$.  Let $$A\pr\xi\lambda\Odef\bigcup
\{\bar {s}_m(\xi ):\,m>l\}.$$

\noindent {\it This completes the construction of} $A\pr\xi\lambda${\it .}

It remains to show that the induction hypotheses ($\aleph_1$~1) through ($\aleph_1$~4) hold
for these $A\pr\xi\lambda$.  Unless otherwise noted, $p_m=\otss mrsn$.

\Claim($\aleph_1$ 1).  For each $\xi <\lambda$, $\langle A\pr\xi\alpha :\xi
<\alpha\leq\lambda\rangle$ is a tower.  
\showClaim It suffices to show for $\alpha$ with $\xi <\alpha <\lambda$ that 
$A\pr\xi\alpha\alss A\pr\xi\lambda$.  So fix such an $\alpha$.  Since $\prxa A\Sminus 
T_{\lambda}$ is finite, there is an $i_0\in T_{\lambda}$ with $i_0>\sup(\prxa A\Sminus 
T_{\lambda})$.  Let $m>\sup(g^{-1}\{ i_0,\op\xi 0,\op\alpha 1\})$.  This choice ensures that 
$\xi\in dom(s_m)$ and $\alpha\in dom(r_m )$.  We show that $\prxa A\Sminus n_m\subseteq 
A\pr\xi\lambda$.  Suppose $i\in\prxa A\Sminus n_m$.  By the definition of $i_0$ and since $i\in 
\prxa A\Sminus i_0$ we know $i\in T_{\lambda}$.  Next, $m'\Odef\Inv g(i)>m$, since if $m'\leq 
m$ then $i$ would be in some $\bar {s}_m(\xi )$ because of the construction at stage $ m'$ and 
this implies $i\leq n_m$ contrary to the choice of $i$.  Hence $p_{m'}\leq p_m$ and by the 
definition of extension, in particular by (e5), we may conclude $i\in\bar {s}_{m'}(\xi )$.  
\endClaim 

\Claim($\aleph_1$ 2).  $T_{\lambda}$ is the disjoint union of the $ A\pr\xi\lambda$ with $\xi
<\lambda$.  

\showClaim Only elements of $T_{\lambda}$ can be put in any $A\pr\xi\lambda$.
That $T_{\lambda}$ is contained in the union follows from the third case in the
definition of
$p_m$ and the fact that $T_{\lambda}$ is contained in the image of $ g$.  The property of
disjointness follows from condition (\bigP1) in the definition of \bigP.  \endClaim

\Claim($\aleph_1$ 3).  For each $\beta$ in $\varphi_{\lambda}\Sminus (\eta_0+1)$,
$A\pr{\xi_0}\lambda\cap A\pr{\eta_0}\beta\neq\eset$.  

\showClaim Fix $\beta\in\fsubl\Sminus(\eta_0+1)$ and let $m>g^{-1} (\op\beta 1)$.  Then 
$\beta\in dom(r_m)$ and by (\bigP6) we conclude $\sbar{\xi_0}\cap A\pr{\eta_0}\beta\neq\eset$.  
But $\sbar{\xi_0} \subseteq A\pr{\xi_0}\lambda$ and so $A\pr{\xi_0}\lambda\cap 
A\pr{\eta_0}\beta\neq\eset$ as desired.  \endClaim 

\Claim($\aleph_1$ 4).  If $\alpha <\lambda$ then $\{\xi :\,\prxa A\nsubseteq A\pr\xi\lambda
\}$ is finite.  

\showClaim Temporarily denote $\{\xi :\,\prxa A\nsubseteq A\pr\xi \lambda \}$ by $X(\alpha )$.  
Fix $\alpha <\lambda$ and let $m\Odef g^{-1}(\op\alpha 1)+1$.  We show $X(\alpha )\subseteq 
f\Simage(r_m(\alpha )+1)$.  

Recall that $f:\omega\to\lambda$ is a bijection.  It is sufficient to show that
$\lambda\Sminus f\Simage(r_m(\alpha )+1)\cap X(\alpha )=\eset$.  Fix $\xi$ such that
$f^{-1}(\xi )>r_m(\alpha )$ and fix $i\in\prxa A$.  Note by (\bigP2C) that $i\in T_{\lambda}$.
Let $l>\max(g^{-1}(i),g^{-1}\langle\xi ,\,0\rangle ,m)$.

Since $p_l\leq p_m$, we have $r_l(\alpha )=r_m(\alpha )$.  Further, $\xi\in dom(s_l)$ and $
i<n_m$.  By (\bigP2A) we can conclude that $i\in\bar {s}_m(\xi )$, hence $i\in
A\pr\xi\lambda$.  Since $i$ was an arbitrary element of $\prxa A$, we have $\prxa A\subseteq
A\pr\xi\lambda$, so $\xi\notin X(\alpha )$.  This gives the desired result.  \endClaim

This completes the proof that the induction hypotheses are satisfied through the limit stage.
Continuing the construction through the countable ordinals gives the desired
$\omega_1\times\omega_1$ matrix.  \endThm

We are now in a position to apply the theorem to obtain information about the gap cohomology.
According to Lemma \the\LocalLLemma, if H2 holds of the ladder system
$\ladder\varphi \lambda\Lambda$ and we apply the $\aleph_1$ gap construction with this system,
then for each pair $\eta <\xi <\omega_1$, the pair of towers $\langle A(\eta ,\, \alpha
),\,A(\xi ,\,\alpha )\,:\,\xi <\alpha <\omega_1\rangle$ is a Hausdorff gap.  This is where the
appearance of $\eta$ in Lemma \the\LocalLLemma\ is used and is necessary because
of the form of the gaps constructed.

Suppose that $X\subseteq\omega_1$ is non-empty and not equal to $ \omega_1$.  Define $$A\pr
X\alpha\Odef\bigcup \{\prxa A:\,\xi\in X\cap\alpha \}.$$ 

\Claim The tower $A(X)\Odef\langle A\pr X\alpha :\,\alpha\in\omega_ 1\rangle$ is a gap in $T$.  

\showClaim Assume that $\mu_0\in X$ and $\mu_1\notin X$.  In fact, if $A(X)$ is truly a 
subtower, it is easy to see that it is a gap in $T$.  This is simply because the tower 
$A(\mu_0)$ is a subtower of $A(X)$ while $A(\mu_1)$ is a subtower of its levelwise complement.  
Any separation of $A(X)$ from its levelwise complement would provide a separation of $ 
A(\mu_0)$ from $A(\mu_1)$.  

So it remains to show that $A(X)$ is a tower.  Fix $\beta <\alpha$.  There are three things to
check:
\item{(a)} $A\pr X\beta\alss A\pr X\alpha$,
\item{(b)} $A\pr X\beta\notaeq A\pr X\alpha$, (i.e.  $A\pr X\alpha\Sminus A\pr X\beta$ is
infinite) and
\item{(c)} The set $A\pr X\alpha\cap T_{\beta}\Sminus A\pr X\beta$ is finite.

(a) and (b) show $A(X)$ is a tower while (c) ensures $A(X)$ satisfies the ``faithful
restriction'' clause in the definition of subtower.

For (a) we need to check that $A\pr X\beta\Sminus A\pr X\alpha$ is finite.  This set is equal
to $$\bigcup_{\eta\in X\cap\beta}\left(\prhb A\Sminus\bigcup_{\xi\in X\cap\alpha}\prxa
A\right).$$ Note that this displayed set is contained in $\bigcup \{\prxb A\Sminus\prxa
A\,:\,\xi\in X\cap\beta \}$ since we have made the set being subtracted smaller.  For each
$\xi\in X$, we have $\prxb A\alss\prxa A$.  However, if $X\cap\beta$ is infinite, there is the
possibility that the finite sets $\prxb A\Sminus\prxa A$ might accumulate to an infinite set.
This is the reason the construction insists that $\prxb A\subseteq\prxa A$ (not almost
containment) ``most of the time''.  In particular, for this fixed $\alpha$ and $\beta$, for
all but finitely many $\xi$ we have $\prxb A\Sminus\prxa A=\eset$, and hence the displayed set
is finite.  This gives (a).

For (b), we need to show $A\pr X\alpha\Sminus A\pr X\beta$ is infinite.  Since
$A\pr{\mu_0}\alpha\subseteq A\pr X\alpha$, it suffices to show that $A\pr{\mu_0}\alpha\Sminus
T_{\beta}$ is infinite because $A\pr X\beta$ is contained in $T_{\beta}$.  Recall that the set
$\{\eta :\,\prhb A\nsubseteq A\pr\eta\alpha \}$ is finite and that the levels of the matrix
are disjoint families.  These, together with the fact that $T_{\beta}\alss T_{\alpha}$, imply
that $A\pr{\mu_ 0}\alpha\cap T_{\beta}=^{*}A\pr{\mu_0}\beta$.  But $A\pr{\mu_0}\alpha\Sminus
A\pr{\mu_0}\beta$ is an infinite set.

For (c), the set of interest is contained in $\bigcup \{\prxa A\cap T_{\beta}\Sminus\prxb
A\}$.  Recall that for all but finitely many $ \xi$ we have $\prxa A\cap T_{\beta}=\prxb A$,
and, as noted in (b), immediately above, for all $ \xi$ we have $\prxa A\cap
T_{\beta}=^{*}\prxb A$.  Hence only finitely many sets contribute to the noted union, and each
only a finite amount.  Thus $A(X)$ satisfies the ``faithful restriction'' requirement.
\endClaim

This claim leads to the following 
\Cor  It is consistent with the statement $2^{\aleph_1}> 2^{\aleph_0}$ 
that the gap cohomology group have cardinality $2^{\aleph_1}$.  

\showCor We have shown that H2 is consistent in the preceding section.  Build the
$\omega_1\times\omega_1$ matrix in the $\aleph_1$ Gap Theorem with the ladder system in (G2)
satisfying H2.  As noted in the corollary to the $\aleph_0$ Gap Theorem, if $X$ and $Y$ are
non-trivial subsets of $ \omega_1$ then $A(X)\symdiff A(Y)=A(X\symdiff Y)$.  In the discussion
of the gap cohomology, recall that two gaps represent different cohomology classes just in
case their levelwise symmetric difference is again a gap.  By the immediately preceding claim,
this is the case whenever $X\neq Y$ and $X\neq \omega_1\Sminus Y$.  \endCor

This ends the discussion of cohomology and gaps for this article.  In attempting to
settled the issue of, for example, the possible size of the gap cohomology group, it may
be useful to look at definable properties of gaps.  One example of this is tight gaps
which are described in [Sc] and [Ra].  The next section introduces another example,
incollapsible gaps.

\Section{Incollapsible Gaps}
Let $\langle A,\,B\rangle$ be a Hausdorff gap.  We ask:  On what subsets of $
\omega$
does the restriction of $\langle A,B\rangle$ remain a gap?  Note that the question
only makes sense for subsets on which $\langle A,B\rangle$ remains a pregap.
For a tower, $A$, in $\Cal P(\omega )$ let $A\rst Z$ for $Z\subseteq
\omega$ be the tower
$\langle A_{\alpha}\cap Z\,:\,\alpha\in\omega_1\rangle$.
\DDef{Collapses (a gap).} Say that a gap, $\langle A,B\rangle$, {\it collapses\/} on
$Z\subseteq\omega$ if and only if $\langle A\rst Z,B\rst Z\rangle$ is a pregap  but not a gap.
That is, there is some $Y\subseteq Z$ that fills $\langle A\rst Z
,B\rst Z\rangle$.
\DDef{Incollapsible (a gap).} Say a gap, $\langle A,B\rangle$, is {\it incollapsible\/} if
$$(\forall Z\subseteq\omega )\,\langle A,B\rangle\text{ does not collapse on }
Z.$$

\DDef {IG} Let IG be the statement:  There is an incollapsible gap.
\Thm(Incollapsible Gaps.)  IG is independent of ZFC.
\showThm
We show that MA + $\neg$CH $\vdash$ $\neg$IG and CH $\vdash$ IG.  Similar
investigations have been undertaken in [K-vD-vM] who show
under MA + $\neg$CH  that for any gap there is an infinite  proper 
subset of $\omega$ on which the gap remains a gap when restricted.

\noindent {\it Sketch of proof of\/} MA + $\neg$CH $\vdash$ $\neg$IG.

The ``obvious'' partial order works.  Fix $\langle A,B\rangle$, a gap.  Let
$$\Bbb P\Odef\{\langle z,\,y,\,s,\,t\rangle\,:\,(\exists n\in\omega) z,y\in{}^n\{0,1\}\,
\land\,s,t\in [\omega_1]^{<\omega}\}.$$
$z$ will build a set $Z$, while $y$ will build a set $Y\subseteq
Z$.  The set $s$ is a list of ordinals, $\alpha$, for which we promise to keep
$A_\alpha\cap Z$ inside $Y$ ``from now on'' and $t$ is a list of ordinals
keeping
$B_\alpha\cap Z$ out of $Y$ ``from now on''.

The goal is to
have $\langle A\rst Z,B\rst Z\rangle$ a pregap filled by $Y$.  Recall the notation
from previous chapters that $\bar {z}$ is the set $z^{-1}\{1\}$, etc.
Let $p=\langle z,\,y,\,s,\,t\rangle$ and
$p'=\langle z',\,y',\,s',\,t'\rangle$ and define $\leq$ on $\Bbb P$ by
$p\leq p'$ if and only if
\item{(1)} $z\supseteq z'$, $y\supseteq y'$, $s\supseteq s'$ and $
t\supseteq t'$.
\item{(2)} (With the obvious notational conventions:)
$(\forall i\in n\Sminus n')$
\itemitem{(i)} $(\forall\alpha\in s')\,i\in A_{\alpha}\,\land\,i\in
\bar {z}\implies i\in\bar {y}$,
\itemitem{(ii)} $(\forall\alpha\in t')\,i\in B_{\alpha}\,\land\,i
\in\bar {z}\implies i\notin\bar {y}$.
\Claim$\Bbb P$ has the c.c.c.  (In fact, it is $\sigma$-centered.)
\showClaim\  Fix an uncountable subcollection of $\Bbb P$, $S.$
By thinning, we may assume all the $z$ and $y$ components are the
same, independent of the condition chosen from $S$.  The
componentwise union of any pair of these is a common extension.

\endClaim

We now describe the dense sets that ensure $\langle A,B\rangle$ will collapse
as desired.  Let $\beta <\alpha <\omega_1$ and $m\in\omega$ be fixed.  Define
$$D_{\alpha\beta}^m\Odef\{p=\langle z,\,y,\,s,\,t\rangle\,:\,\vert
\bar {z}\cap A_{\alpha}\Sminus A_{\beta}\vert\geq m\}.$$
\Claim$D^m_{\alpha\beta}$ is dense for all $\alpha$, $\beta$ and $
m$.
\showClaim\  By induction on $m.$ The claim is clear for $m=0$.
Suppose it is true for $m$.  Fix $p_0\in \Bbb P$ and $p_1\leq p_0$ (with the
obvious notational extensions) so that
$\vert\bar {z}_1\cap A_{\alpha}\Sminus A_{\beta}\vert\geq m$.  Let $
r\in\omega$ be so large that
$$\left(\bigcup_{\alpha'\in s_1\cup t_1}A_{\alpha'}\Sminus r\right
)\cap\left(\bigcup_{\alpha'\in s_1\cup t_1}B_{\alpha'}\Sminus r\right
)=\emptyset$$
which is possible since $s_1\cup t_1$ is a finite set.  Fix $n\in A_{\alpha}\Sminus 
A_{\beta}$, an infinite set, such that $n>r,n_1$.  Let $ z\supseteq z_1$ and $n\in\bar 
{z}$, and let $y\supseteq y_1$ with $n\in\bar {y}$.  Then $\langle z,\, 
y,\,s,\,t\rangle\leq p_0$ and $\vert\bar {z}\cap A_{\alpha}\Sminus A_{\beta}\vert\geq 
m+1$.  \endClaim 

The point of this claim is that $A\rst Z$
will still be an $\omega_1$ tower.
A similar argument shows the same for $B\rst Z$ and so
$\langle A\rst Z,B\rst Z\rangle$ remains a pregap.  To show that
the set $Y$ separates the restricted pregap, note that if
$\alpha\in s$ and $\langle z,\,y,\,s,\,t\rangle\in G$ for $G$ generic then $
Y\supseteq^{*}(Z\cap A_{\alpha})$.  So it
suffices to show the set $E_{\alpha}\Odef\{\langle z,\,y,\,s,\,t\rangle
\in \Bbb P\,:\,\alpha\in s\}$ is dense which is left to the reader.

\noindent {\it Sketch of proof of\/} CH $\vdash$ IG.

This is a standard diagonalization argument on $[\omega ]^{\aleph_
0}$ (the infinite subsets of $\omega$) which under
CH has cardinality $\aleph_1$.  Fix an enumeration of $[\omega ]^{
\aleph_0}$, $Y_{\alpha}$ for $\alpha <\omega_1$.
Recursively
define $A_{\alpha}$ and $B_{\alpha}$ so that $Y_{\alpha}$ will not fill $
\langle A\rst Z,B\rst Z\rangle$ for
any $Z\supseteq Y_{\alpha}$.
Suppose $A_{\beta}$ and $B_{\beta}$ are defined for $\beta <\alpha$ and form a pregap.  If
for some $\beta <\alpha$ we have $Y_{\alpha}\subseteq^{*}A_{\beta}$ then there is nothing to worry
about.  Otherwise $Y_{\alpha}\Sminus A_{\beta}$ is infinite for each $
\beta <\alpha$.  As the set of $A_{\beta}$'s
is an increasing chain, the collection of $Y_{\alpha}\Sminus A_{\beta}$ has the
strong finite intersection property.  So there is an infinite set
$X\subseteq^{*}Y_{\alpha}\Sminus A_{\beta}$ and thus $X$ is almost disjoint from each $
A_{\beta}$.

Since $\alpha$ is countable, there is a set $X'$ such that for all $
\beta <\alpha$,
$X'\supsetneq^{*}B_{\beta}$ and $X'\cap A_{\beta}=^{*}\emptyset$.  Let $
B_{\alpha}\Odef X\cup X'$ and notice
now that $Y_{\alpha}$ cannot separate $A$ from $B$ since $B_{\alpha}
\cap Y_{\alpha}$ is infinite.
Let $A_{\alpha}$ be any set such
that $A_{\alpha}\supsetneq^{*}A_{\beta}$ for $\beta <\alpha$ and $
A_{\alpha}$ is almost disjoint from
$B_{\alpha}$.
\endThm


\vskip.3truein plus .2 truein minus .1 truein
\centerline{\bf References}
\vskip .1truein plus .1truein minus .05truein
\parskip=\MediumSpaceBase plus 2pt minus 1pt
\baselineskip=\MediumSpaceBase

[Ba]    J.~Baumgartner,  Applications of the Proper Forcing
Axiom, in: K.~Kunen and 
J.~Vaughan, eds.,   Handbook of set-theoretic topology (North-Holland, New York,   1984).

[Bl]  A.~Blass,   Cohomology detects failures of the Axiom of
Choice,  Trans.\ Amer.\ Math.\ Soc.  279, no.~1 (1983)
257--269.

[De] K.~Devlin, Variations on $\diamondsuit$,  J.\ Symbolic Logic,
 44, no.~1 (1979) 51--58.

[D-S] K. Devlin and S. Shelah, A weak version of $\diamondsuit$ which follows
from $2^{\aleph_0}<2^{\aleph_1}$, Israel J.~Math.~29 (1978) 239--247.

[D-S-V] A.~Dow,   P.~Simon, and   J.~Vaughan,  Strong
Homology and the Proper Forcing Axiom,   Proc.\ Amer.\
Math.\ Soc., 106, no.~3 (1989),    821--828.

[Dr] F.~Drake, Set theory:  An introduction to large cardinals
(North-Holland, Amsterdam-New York, 1974).

[E-M] P.~Eklof and A.~Mekler, Almost free modules:  Set
theoretic methods (North-Holland, Amsterdam-New York, 1990).

%

[F]  D.~Fremlin, Consequences of Martin's Axiom  (Cambridge
University Press, Cambridge, 1984).


[Ha] F.~Hausdorff, Summen von $\aleph_1$ Mengen, Fund.\ Math., 
26 (1936) 241--255.

%

[H-Y] J.~Hocking and G.~Young, Topology  (Dover, New York, 1961).

%

[Jn] C.~Jensen, Les foncteurs d\'eriv\'es de $\varprojlim$
 et leurs applications en th\'eorie des modules, Lecture Notes in
Mathematics, v.~254 (Springer-Verlag, New York, 1972).

[Ku] K.~Kunen, Set theory: An introduction to independence proofs
(Elsevier North-Holland, Amsterdam, 1980).

%

%

%

%

%

[Os] A.~Ostaszewski, On countably compact, perfectly normal spaces,
J.\ London Math.\ Soc., 2nd series, 14, no.~4 (1976) 505--516.

[Ra] M.~Rabus, Tight gaps in $\Cal P(\omega)$, preprint.

%

[Sc] M.~Scheepers, Gaps in $({}^\omega\omega, \prec)$,  Israel Mathematical Conference
Proceedings,  6 (1993) 439--561.

%

%

[T1]  S.~Todor\v cevi\'c, Partition problems in topology
(Contemp.\ Mathematics, 84, Amer.\ Math.\ Soc.,
Providence, RI, 1989).

[T2]  S.~Todor\v cevi\'c, Remarks on MA and CH, Can.\ J.\ Math.,
43 (1991) 832--841.

%

%

[Vi] J.~Vick,   Homology theory, (Academic Press, New York, 1973).

\end